\documentclass[12pt,reqno]{amsart}

\usepackage{AKstyle}

\numberwithin{equation}{section}

\usepackage[final, 
color]{showkeys}
\definecolor{refkey}{gray}{.85}
\definecolor{labelkey}{gray}{.85}

\usepackage[pagebackref=true, colorlinks]{hyperref}

\hypersetup{pdffitwindow=true,linkcolor=blue,citecolor=blue,urlcolor=blue,menucolor=blue}

\usepackage{comment}

\newcommand\Bpow{\fb}

\newcommand\TwenOvQpow{\frac34\frac\Bpow R} 

\newcommand\Ais{50}
\newcommand\Aispo{51}
\newcommand\gdIs{5/312} 
\newcommand\ogdIs{312/5} 
\newcommand\Aapprox{37.93} 
\newcommand\OmDelApproxx{0.984} 
\newcommand\OmDelIs{307/312} 

\newcommand\gdquart{{13/5}}
\newcommand\gdhalf{{26/5}}
\newcommand\gdone{{52/5}}
\newcommand\gdtwo{{104/5}}
\newcommand\gdover{{5/52}}

\newcommand\gdofour{{5/208}}

\newcommand\gba{\frac1{20}}
\newcommand\gbb{\frac1{10}}

\newcommand\gbaa{\frac1{25}}
\newcommand\gbbb{\frac1{5}}

\newcommand\cIm{\frac1{100}}
\newcommand\cIp{\frac{99}{100}}

\begin{document}

\author{Jean Bourgain}
\thanks{Bourgain is partially supported by NSF grant DMS-0808042.}
\email{bourgain@ias.edu}
\address{IAS, Princeton, NJ}
\author{Alex Kontorovich}
\thanks{Kontorovich is partially supported by  NSF grants DMS-1209373, DMS-1064214 and DMS-
1001252.}
\email{alex.kontorovich@yale.edu}
\address{Stony Brook University, Stony Brook, NY}
\curraddr{Yale University, New Haven, CT}

\title{On 
Zaremba's conjecture
}

\begin{abstract}
Zaremba's 1971 conjecture predicts that every integer appears as the
denominator
 of a finite continued fraction whose
  partial quotients 
  are 
  bounded
  by an absolute constant.
   We confirm this conjecture for a set of density one. 
\end{abstract}
\date{\today}
\maketitle
\tableofcontents


\section{Introduction}

\subsection{Statements of the Main Theorems}

For a fixed finite set $\cA\subset\N$, which we call an {\bf alphabet},
let $\fC_{\cA}$ denote the collection of all $x\in(0,1)$ whose continued fraction expansion
$$
x
=[a_{1},a_{2},\dots,a_{k},\dots]
=\cfrac{1}{a_{1}+\cfrac{1}{a_{2}+\ddots+\cfrac{1}{a_{k}+\ddots}}}
,
$$
has all partial quotients $a_{j}$ belonging to the alphabet $\cA$. Any $x\in\fC_{\cA}$ is uniformly badly approximable, in the sense that its partial quotients $a_{j}$ are all bounded by
$$
A:=\max\cA.
$$
When $A$ is an absolute constant, we call such a number {\bf absolutely 
Diophantine} (of height $A$). That is, when we speak of numbers being absolutely Diophantine, the height $A$ is fixed in advance.

Let $\fR_{\cA}$ 
denote the set of 
partial convergence to $\fC_{\cA}$, that is,
$$
\fR_{\cA}:=\left\{\frac bd=[a_{1},a_{2},\dots,a_{k}]:0<b<d,\ (b,d)=1,\text{ and } \forall j,\ a_{j}\in\cA\right\}
,
$$ 
and let $\fD_{\cA}\subset\N$ be the set of denominators of fractions in $\fR_{\cA}$,
$$
\fD_{\cA}:=\left\{d\in\N:\exists(b,d)=1\text{ with }\frac bd\in\fR_{\cA} 
\right\}
.
$$

In 1971, S. K. Zaremba formulated the following assertion.

\begin{conj}[{Zaremba \cite[p. 76]{Zaremba1972}}]\label{conj:Z}
Every positive integer is the denominator of a reduced absolutely 
Diophantine
fraction.
\end{conj}
That is, the conjecture predicts the existence of some integer $A>1$ so that 
$$
\fD_{\{1,2,\dots, A\}}=\N.
$$

Zaremba's conjecture has important applications to numerical integration and pseudorandom number generation, producing 
collections of points 
of optimal discrepancy; see e.g. the surveys \cite{Niederreiter1978} and \cite{Kontorovich2013}.
Our main result is the following
\begin{thm}\label{thm:main}
Almost every positive integer is the denominator of a reduced 
absolutely 
Diophantine
fraction. That is, there exists an effectively computable $A>1$ so that 
$$
\frac1N\#(\fD_{\{1,2,\dots,A\}}\cap[1,N])\to1,
$$
as $N\to\infty$.
\end{thm}

A more refined conjecture was stated by Hensley in 1996. The set $\fC_{\cA}\subset(0,1)$ is a Cantor-like fractal; let 
$$
\gd_{\cA}:=\text{H.dim}(\fC_{\cA})\ \in \ [0,1]
$$ 
be its Hausdorff dimension. This dimension can be $0$ only if  $|\cA|=1$; since we assume $\cA$ is finite, $\gd_{\cA}<1$.
Allowing a finite number of exceptions in Zaremba's conjecture, Hensley asserts the following.

\begin{conj}[{Hensley \cite[Conjecture 3, p.16]{Hensley1996}}]\label{conj:Hensley}
The set of denominators 
$
\fD_{\cA}
$
contains every sufficiently large integer if and only if the
corresponding dimension $\gd_{\cA}$ 
exceeds $1/2$.
\end{conj}

As stated, Hensley's conjecture is false.
For example, consider the 
alphabet $\cA=\{2,4,6,8,10\}$.
By implementing an algorithm due to Jenkinson and Pollicott \cite{JenkinsonPollicott2001},%
\footnote{
The program is available 
at  
\url{http://math.sunysb.edu/\~alexk/HausdorffZaremba.nb}.
}
we have 
estimated
its dimension to be
$
\gd_{\cA}\approx 0.517>1/2.
$
Nevertheless,
arbitrarily large numbers are
  missing from $\fD_{\cA}$.
Indeed, it is elementary to verify that
\be\label{eq:fDcAmod4}
\fD_{\cA}(\mod 4)\equiv\{0,1,2\},
\ee
see Remark \ref{rmk:StongApp}. 

We propose the following alternative to Hensley's conjecture, borrowing language from Hilbert's $11^{th}$ problem on representations of numbers by quadratic forms. 
We call an integer $d$ 
 {\bf admissible} (for $\cA$) if it 
 passes all
 finite
  local obstructions:
\be\label{eq:AdmissDef}
\forall q>1,\
d\in\fD_{\cA}(\mod q)
.
\ee
\begin{rmk}\label{rmk:q}
Admissibility can be checked using only one modulus $q=q(\cA)$, see Remark \ref{rmk:StongApp}.
\end{rmk}
Let  
$
\fA_{\cA}
$ 
denote  the set of all admissible numbers,
$$
\fA_{\cA}:=\{d\in\Z: \eqref{eq:AdmissDef}\text{ holds}\}
.
$$
We say $d$ is {\bf represented} (by $\cA$) if $d\in\fD_{\cA}$. The {\bf multiplicity} of a denominator $d$ is the number of coprime numerators $0<b<d$ with $b/d\in\fR_{\cA}$. Clearly $d$ is represented if and only if its multiplicity is positive.

\begin{conj}\label{conj:BK}
If
the dimension $\gd_{\cA}$ exceeds $1/2$,
then
the set of denominators 
$
\fD_{\cA}
$
contains every sufficiently large admissible integer.
\end{conj}

We interpret this conjecture as a local-global principle, where the dimension condition and ``sufficiently large'' are local obstructions at infinity.
Theorem \ref{thm:main} follows from the following
more refined
 approximation to Conjecture \ref{conj:BK}.
\begin{thm}\label{thm:mainN}
There exists an effectively computable constant $\gd_{0}<1$ so that
if
the dimension
$\gd_{\cA}$ exceeds $\gd_{0}$,
then
the set of denominators 
$
\fD_{\cA}
$
contains almost every admissible integer.
%
More precisely, there is a constant $c=c(\cA)>0$ so that 
\be\label{eq:smExcept}
{\#(\fD_{\cA}\cap[N/2,N])
\over
\#(\fA_{\cA}\cap[N/2,N])
}
=
1+O\left(N^{-c/\log\log N}\right)
,
\ee
as $N\to\infty$.
Furthermore, each $d$ produced above appears with multiplicity 
\be\label{eq:mult}
\gg N^{2\gd_{\cA}-\frac{1001}{1000}}
.
\ee
The constants $\gd_{0}$ and $c$ are effectively computable, and the implied constants above 
depend only on $\cA$.
\end{thm}

Some remarks are in order.

\begin{rmk}
There exist alphabets $\cA$ with $\gd_{\cA}$ arbitrarily close to $1$, so Theorem \ref{thm:mainN} is not vacuous. Indeed, Hensley \cite{Hensley1992} gives the asymptotic expansion
\be\label{eq:delHens}
\gd_{\{1,2,\dots,A\}}=1-{6\over \pi^{2}A}-{72\log A\over \pi^{4}A^{2}}+O\left(\frac1{A^{2}}\right).
\ee
\end{rmk}

\begin{rmk}
The 
number
$1/2$ in Conjectures \ref{conj:Hensley}  and  \ref{conj:BK}  
cannot be reduced. Hensley \cite{Hensley1989} showed that the truncated set of rationals
$$
\fR_{\cA}(N):=\left\{\frac bd\in\fR_{\cA}:(b,d)=1,\ 0<b<d<N\right\}
$$
has cardinality
\be\label{eq:2dim}
\#\fR_{\cA}(N)\
\asymp
\ N^{2\gd_{\cA}}
,
\ee
whence
it follows immediately that
$$
\#(\fD_{\cA}\cap[1,N])\ll N^{2\gd_{\cA}}.
$$
Thus if $\gd_{\cA}<1/2$, then 
 certainly
  $\fD_{\cA}$ is too thin a subset of the integers to contain even one admissible arithmetic progression. 
\end{rmk}

\begin{rmk}
The best previously known estimate 
\be\label{eq:DANpower}
\#(\fD_{\cA}\cap[1,N]) \gg N^{\gd_{\cA}}.
\ee
was proved by Hensley  \cite[Theorem 3.2]{Hensley2006}, and follows easily from his estimate \eqref{eq:2dim}. In particular, as long as $|\cA|>1$, the set $\fD_{\cA}$ grows at least polynomially.
Moreover, taking $\cA$ large so that $\gd_{\cA}>1-\vep$, one 
can already
produce 
at least $N^{1-\vep}$ denominators in $\fD_{\cA}$ up to $N$.
\end{rmk}

\begin{rmk}\label{rmk:allQ}
We
explain
in
Remark \ref{rmk:StongApp} below
that
for any $A\ge2$, 
 the alphabet $\{1,2,\dots, A\}$ has no finite local obstructions, that is, $\fA_{\{1,2,\dots,A\}}=\Z$. This 
 is why
  the statement of Theorem \ref{thm:main}
  needs no mention of admissibility.
Moreover the dimension $\gd_{\{1,2\}}$ is known \cite{Good1941, Bumby1985, JenkinsonPollicott2001} to be
\be\label{eq:del2Is}
\gd_{\{1,2\}}\approx0.531\cdots
,
\ee
which obviously exceeds $1/2$. 
Conjecture \ref{conj:BK} 
then implies that $\fD_{\{1,2\}}$ already contains every sufficiently large number, as was conjectured by Hensley \cite{Hensley1996}.
\end{rmk}

\begin{rmk}
An earlier version\footnote{\url{http://arxiv.org/abs/1107.3776v1}} of this paper also proved two weaker results 
made obsolete
by Theorem \ref{thm:main}, namely that for sufficiently large $\gd_{\cA}$, $(i)$ $\fD_{\cA}$ contains a positive proportion of numbers, and $(ii)$ that $\fD_{\cA}$ contains almost every admissible number, without giving the rate in \eqref{eq:smExcept}.
At the request of the referee to shorten the paper, we have removed these intermediary results (and of course the methods used to obtain them). We invite the interested reader to peruse the original arxiv posting for the details.
Note also that some results of this paper have been announced in \cite{BourgainKontorovich2011}.
\end{rmk}

\begin{rmk}
The value of $\gd_{0}$ in Theorem \ref{thm:main} coming from our proof is
\be\label{eq:gdAllow}
\gd_{0}=\OmDelIs\approx\OmDelApproxx.
\ee
We have made no effort to optimize this quantity, as can surely be done with a modicum of effort. In fact, Frolenkov and Kan\footnote{\url{http://arxiv.org/abs/1303.3968v1}} have since
sharpened our method to prove 
the weaker statement
that $\fD_{\cA}$ contains a positive proportion of numbers whenever $\gd_{\cA}>\gd_{0}$ with
the 
improved
range
 $\gd_{0}=5/6\approx0.833$. It does not seem likely that our methods can achieve the full range $\gd_{0}=1/2$ without significant new ideas. 
We have estimated the dimension $\gd_{\cA}$ corresponding to 
 the alphabet $\cA=\{1,2,\dots,49,50\}$
 to be about $0.986$,
 exceeding \eqref{eq:gdAllow}, whereas the alphabet $\cA=\{1,2,3,4,5\}$ is known \cite{Jenkinson2004} to have dimension $\gd_{\cA}>5/6$. 
\end{rmk}

Although
our main result requires
  large dimension, 
 we are 
 also
 able
  to sharpen 
  the best previously known estimate \eqref{eq:DANpower} in the full range $\gd_{\cA}>1/2$.
\begin{thm}\label{thm:gdPlus}
Write $\gd$ for $\gd_{\cA}$. Then for any $\vep>0$,
\be\label{eq:gdPlus}
\#(\fD_{\cA}\cap[1,N])\gg_{\vep} N^{\gd+\frac{(2\gd-1)(1-\gd)}{ 5-\gd}-\vep}
,
\ee
as $N\to\infty$. This bound improves on \eqref{eq:DANpower}, as long as $\gd>1/2$.
\end{thm}

\begin{rmk}
The improvement here is quite modest: for $\cA=\{1,2\}$, the exponent $
\gd_{\cA}\approx 0.531$ in \eqref{eq:del2Is} and \eqref{eq:DANpower} is replaced in \eqref{eq:gdPlus} by $0.537$. We have again made no attempt to optimize the exponent in \eqref{eq:gdPlus}, seeking just any power gain.
\end{rmk}

We state the multiplicity bound \eqref{eq:mult}   to give another application
to pseudorandom numbers. 
Specifically,
in the (homogeneous) linear congruential method, 
optimal conditions require
 a prime $d$ and a primitive root $b(\mod d)$ so that the fraction $b/d$ is absolutely Diophantine (see \cite{Kontorovich2013}). Then the pseudorandom map with modulus $d$ and multiplier $b$, that is, $x\mapsto bx(\mod d)$, has asymptotically optimal serial correlation of pairs.

\begin{thm}\label{thm:primes}
There exist 
infinitely many
primes $d$ with primitive roots $b(\mod d)$ so that the fractions $b/d$ are absolutely Diophantine. 
\end{thm}

The number of such prime $d$ up to $N$ provided by our proof is $\gg N (\log N)^{-2}$.
Theorem \ref{thm:primes} is an easy corollary of Theorem \ref{thm:mainN}.
In fact, if $\cA=\{1,2,\dots,A\}$ has dimension $\gd_{\cA}$ exceeding $\gd_{0}$ as in Theorem \ref{thm:mainN}, then the fractions $b/d$ produced in Theorem \ref{thm:primes} can be taken to have all partial quotients bounded by $A+1$.

\subsection{Reformulation and Admissibility}

It is an old and trivial (but for our purposes crucial) observation that
$$
\frac bd 
=[a_{1},\dots,a_{k}]
$$
is equivalent to
\be\label{eq:1.1}
\mattwo *b*d
=
\mattwo 011{a_{1}}
\mattwo 011{a_{2}}
\cdots
\mattwo 011{a_{k}}
.
\ee
%
This observation will 
allow us to explain all local obstructions, as follows.
In light of \eqref{eq:1.1}, let
$$
\cG_{\cA}\subset\GL(2,\Z)
$$ 
be the semigroup generated 
by the matrices
\be\label{eq:cGgens}
\mattwo 011a
\ee
for $a\in \cA$.
Then the orbit 
\be\label{eq:cOIs}
\cO_{\cA}:=\cG_{\cA}\cdot e_{2}
\ee
of $e_{2}=(0,1)^{t}$ under $\cG_{\cA}$ 
corresponds to 
$\fR_{\cA}$, that is, if 
$
\g=
\bigl( \begin{smallmatrix}
a&b\\ c&d
\end{smallmatrix} \bigr)
,
$
then $\g\cdot e_{2}=(b,d)^{t}.$ 
Moreover, taking the inner product of this orbit with $e_{2}$ picks off the value of $d$, that is $\<\g\cdot e_{2},e_{2}\>=d$, and 
\be\label{eq:Ge2e2}
\<\cO_{\cA},e_{2}\>=\<\cG_{\cA}\cdot e_{2},e_{2}\>
\ee
is precisely $\fD_{\cA}$ (with 
multiplicity).
Zaremba's conjecture can then be reformulated as: For some finite alphabet $\cA$,
$$
\N\subset\<\cG_{\cA}\cdot e_{2},e_{2}\>
.
$$
For convenience we pass from $\cG_{\cA}$ to its determinant one subsemigroup
$$
\G_{\cA}=\cG_{\cA}\cap\SL_{2}\subset\SL_{2}(\Z)
,
$$
which is (freely and finitely) generated by the matrix products 
\be\label{eq:Ggens}
\bp
0&1\\ 1&a
\ep
\cdot
\bp
0&1\\ 1&a'
\ep
,
\ee
for $a,a'\in\cA$.
The orbit $\cO_{\cA}$ is recovered as a finite union of ``coset'' orbits
$$
\cO_{\cA}
=
\G_{\cA}\cdot e_{2}\ 
\cup\ 
\bigcup_{a\in\cA}
\G_{\cA}\cdot\mattwo011a e_{2}
.
$$

\begin{rmk}\label{rmk:StongApp}
It now follows from
Strong Approximation \cite{MatthewsVasersteinWeisfeiler1984} and Goursat's Lemma
(see the discussion in \cite[\S2.2]{Kontorovich2013})
that the reduction of $\G_{\cA}$ mod $q$ is all of $\SL_{2}(q)$, for all $q$
coprime to a certain ``bad'' modulus $\fB$.
Here $\fB$ is effectively computable and depends only on $\G_{\cA}$, that is, on $\cA$. 
Moreover $\fB$ can be chosen so that for any $q\equiv0(\fB)$, the reduction $\G_{\cA}(\mod q)$ is the full pre-image of $\G_{\cA}(\mod \fB)$ under the projection map $\Z/q\to\Z/\fB$. 
From the mod $\fB$ reductions of $\G_{\cA}$, it is elementary to read off the reductions of $\cO_{\cA}$, and hence all finite local obstructions in $\fD_{\cA}$; see Remark \ref{rmk:q}. Moreover, for the alphabet $\cA=\{1,2\}$, it is easy to see that $\fB=1$, that is, $\G_{\cA}(\mod q)$ is already all of $\SL_{2}(q)$, for all $q>1$; see Remark \ref{rmk:allQ}.
Indeed, the {\it group} generated by $\G_{\cA}$ (that is, allowing inverses) is all of $\SL_{2}(\Z)$, and the two have the same projections mod $q$.
Finally, we note
that this is
precisely
 the phenomenon responsible for \eqref{eq:fDcAmod4} in the failure of Hensley's 
 Conjecture \ref{conj:Hensley}.
\end{rmk}

\subsection{An Overview of the Key Ideas}\label{sec:outline}

An observation which we had made in a 
slightly different context
\cite{BourgainKontorovich2010} is that there is a certain bilinear (in fact multilinear) structure to \eqref{eq:Ge2e2}, making the problem amenable to the Hardy-Littlewood circle method via Vinogradov's techniques for estimating bilinear forms. We now outline the key steps.

From now on,  we treat $\cA$ as fixed, dropping it from subscripts, writing $\gd=\gd_{\cA}$, $\G=\G_{\cA}$, etc.
%
In light of \eqref{eq:Ge2e2}, we would like to study the exponential sum 
\be\label{eq:SNbad}
S_{N}(\gt):=\sum_{\g\in\G\atop\|\g\|<N}e(\gt\<\g e_{2},e_{2}\>)
,
\ee
where $\gt\in[0,1]$
and $\|\cdot\|$ is the Frobenius matrix norm, $\|\bigl( \begin{smallmatrix}
a&b\\ c&d
\end{smallmatrix} \bigr)
\|^{2}=a^{2}+b^{2}+c^{2}+d^{2}$.
Then the Fourier coefficient
\be\label{eq:RNd}
R_{N}(d):=\hat S_{N}(d)=\int_{0}^{1}S_{N}(\gt)e(-d\gt)d\gt
=
\sum_{\g\in\G\atop\|\g\|<N}\bo_{\{\<\g e_{2},e_{2}\>=d\}}
\ee
is just 
the
``representation number'' of $d$ up to $N$, that is, its
multiplicity. Of course if $R_{N}(d)>0$, then
$d\in\fD
$.

Note that  by \eqref{eq:2dim},
\be\label{eq:SN0}
S_{N}(0)
=
\sum_{d}R_{N}(d)
=
\sum_{\g\in\G\atop\|\g\|<N}1\
\asymp\ N^{2\gd}
,
\ee
so if almost every $d\in[N/2,N]$ is to be represented without much bias, it 
should
occur  with multiplicity roughly $N^{2\gd-1}$.

Following the circle method, we decompose the integral in \eqref{eq:RNd} into ``major arcs''  and ``minor arcs'', the former referring to modes $\gt$ quite near rationals with small denominators and the latter being the rest:
$$
R_{N}(d)
=
\left(
\int_{\fM}
+\int_{[0,1]\setminus\fM}
\right)
S_{N}(\gt)e(-d\gt)d\gt=\cM_{N}(d)+\cE_{N}(d).
$$
Here $\cM_{N}$ is thought of as a ``main'' term and $\cE_{N}$ is an ``error'' term, and the major arcs $\fM=\fM_{\cQ}$ are given by
\be\label{eq:MajArcsDef}
\fM_{\cQ}=\bigcup_{q<\cQ}\bigcup_{(a,q)=1}\left[\frac aq-\frac {\cQ}N,\frac aq+\frac {\cQ}N\right]
,
\ee
where $\cQ$ is roughly of size $N^{c/\log\log N}$. 

A key ingredient is to show that along the major arcs, $\gt=\frac aq+\gb\in\fM$, 
the function $S_{N}$ essentially splits into two pieces,
\be\label{eq:SNvqvpgb}
S_{N}\left(\frac aq+\gb\right)
\sim
\nu_{q}(a)
\cdot
\vp(\gb)
.
\ee
Here $\nu_{q}$ is a purely modular term and $\vp$ is an archimedean one, which has the right order of magnitude on balls of certain size.
It then follows that the main term $\cM_{N}(d)$ also splits as a ``singular series'' $\fS$ times a ``singular integral'' $\Pi$,
$$
\cM_{N}(d)\sim\fS(d)\Pi_{N}(d),
$$
where $\Pi$ 
gives 
 the expected 
 archimedean contribution,  
 roughly
$$
\Pi_{N}(d)\gg { N^{2\gd-1}}
$$ 
for $d\asymp N$,
and the singular series $\fS$ controls the local obstructions. In particular, if $d\not\in\fA$ is not admissible, then $\fS(d)=0$; otherwise, we have roughly that
$$
\fS(d)\asymp
\prod_{p\nmid d}
\left(
1+
{1\over p^{2}-1}
\right)
\prod_{p\mid d}
\left(
1-
{1\over p+1}
\right)
\gg
{1\over\log\log d}
.
$$

The main ingredient in
 proving 
\eqref{eq:SNvqvpgb} 
is
 the renewal method in the  thermodynamic formalism
 of Ruelle transfer operators
  (see Lalley \cite{Lalley1989}), 
  and
the extension to    ``congruence'' 
such
established by Bourgain-Gamburd-Sarnak in 
\cite{BourgainGamburdSarnak2011}. (We need here not just square-free but
arbitrary moduli $q$, and must also use the work of Bourgain-Varju \cite{BourgainVarju2011}.)


With the major arcs controlled, if we could prove that the errors are individually bounded, $|\cE_{N}(d)|\ll N^{2\gd-1-\vep}$, say, then we would conclude the full Conjecture \ref{conj:BK}. We are not able to establish control of this quality individually, but do succeed on average, proving essentially that

\be\label{eq:SNN4gd}
\sum_{d\asymp N}|\cE_{N}(d)|^{2}
\ll {N^{4\gd-1-c/\log\log N}}
,
\ee
from which Theorem \ref{thm:mainN} follows by a standard argument.

Bounds of this type will follow from bounds on
\be\label{eq:SNonWQK}
\int_{W_{Q,K}}|S_{N}(\gt)|^{2}d\gt,
\ee
where
we have decomposed the minor arcs $[0,1]\setminus \fM$ into the dyadic regions 
\be\label{eq:WQKis}
W_{Q,K}:=\left\{
\gt=\frac aq+\gb :
 q\asymp Q, (a,q)=1,
  |\gb| \asymp \frac KN
\right\}
.
\ee
By Dirichlet's approximation theorem, the parameters $Q$ and $K$ vary in the range
 $Q<N^{1/2}$ and $K<{N^{1/2}\over Q}$.

Unfortunately, we do not know how to obtain such strong bounds for the function $S_{N}$ as defined in \eqref{eq:SNbad}. 
But taking a cue from Vinogradov (as we did in \cite{BourgainKontorovich2010}), we work with a different function:
\be\label{eq:SNno2}
S_{N}(\gt)=
\sum_{\g_{1}\in\G\atop\|\g_{1}\|\asymp N^{1/2}}
\sum_{\g_{2}\in\G\atop\|\g_{2}\|\asymp N^{1/2}}
e(\gt\<\g_{1}\g_{2} e_{2},e_{2}\>)
,
\ee
say. Since $\G$ is a semigroup, this modified function, or rather its Fourier transform,
continues to 
 capture elements of $\fD$. Moreover, the bilinear nature of the problem, namely that 
$\<\g_{1}\g_{2} e_{2},e_{2}\>=\<\g_{2} e_{2},{}^{t}\g_{1} e_{2}\>$, allows us
to separate variables.

It is here in the separation of variables and application of Cauchy-Schwarz that we replace the thin semigroup $\G$ with all of $\SL_{2}(\Z)$, a loss 
we can only tolerate if the dimension $\gd$ is large, at least some $\gd_{0}$. 
We are then lead to a more classical
 setting,
and
 in certain large ranges of the pair $(Q,K)$ in \eqref{eq:WQKis}, we can obtain the 
requisite
 cancelation.
  For slightly smaller values of $(Q,K)$, it is beneficial to decompose the sum further as
$$
S_{N}(\gt)=
\sum_{\g_{1}\in\G\atop\|\g_{1}\|\asymp N^{1/2}}
\sum_{\g_{2}\in\G\atop\|\g_{2}\|\asymp N^{1/4}}
\sum_{\g_{3}\in\G\atop\|\g_{3}\|\asymp N^{1/4}}
e(\gt\<\g_{1}\g_{2}\g_{3} e_{2},e_{2}\>)
.
$$
Continuing in this way, we
handle every 
conceivable 
range of $(Q,K)$ by considering a sum of the form
\be\label{eq:SG1GJ}
S_{N}(\gt)=
\sum_{\g_{1}\in\G\atop\|\g_{1}\|\asymp N^{1/2}}
\sum_{\g_{2}\in\G\atop\|\g_{2}\|\asymp N^{1/4}}
\cdots
\sum_{\g_{J}\in\G\atop\|\g_{J}\|\asymp N^{1/2^{J}}}
e(\gt\<\g_{1}\g_{2}\cdots\g_{J} e_{2},e_{2}\>)
,
\ee
where $J\asymp \log\log N$, so that $\g_{J}$ is of large but constant size (independent of $N$).\footnote{We could take $J$ even a bit smaller, but choose not to for the sake of exposition.}
Unfortunately, another problem has crept up: we can no longer control the size of the long product $\g_{1}\cdots\g_{J}$, which could have norm as large as $N(\log N)^{C}$.

To remedy this situation, we develop a bit of elementary linear algebra for $\G$, showing that if the expanding vectors of two matrices are close, then their eigenvalues behave nearly multiplicatively. This forces us to concoct,
for each $j=1,\dots, J$, a
certain special subset $\Xi_{j}\subset \{\g\in\G:\|\g\|\asymp N^{1/2^{j}}\}$, all the elements of which have expanding eigenvectors 
pointing 
near a common direction (independent of $j$). 
We then simply use the pigeonhole principle to make sure all elements of $\Xi_{j}$ have almost the same eigenvalues.
Moreover, we need to ensure that
the representation 
$\g
=
\g_{1}\g_{2}\cdots\g_{J}
$ 
in \eqref{eq:SG1GJ}
is
unique that is, if
$$
\g_{1}\g_{2}\cdots\g_{J}
=
\g_{1}'\g_{2}'\cdots\g_{J}'
$$
with $\g_{j},\g_{j}'\in\Xi_{j}$, then $\g_{j}=\g_{j}'$ for all $j$.
We do this by forcing
 each $\g_{j}\in\Xi_{j}$ to have the same size in the wordlength metric, again by pigeonhole.  

Then the large product ensemble
$$
\Xi_{1}\cdot\Xi_{2}\cdots\Xi_{J}
,
$$
is a 
good substitute for $\{\g\in\G:\|\g\|\asymp N\}$
to handle the minor arcs. 
Unfortunately, the concocted sets $\Xi_{j}$ are no longer amenable to the major arc methods!
We 
rectify this by constructing
a certain tiny  set $\aleph\subset\G$ with good modular/archimedean distribution properties, 
and 
prepending it to
the product, forming
\be\label{eq:gWNtilde}
\gW_{N}
=
\aleph\ \tilde\Xi_{1} \Xi_{2}\cdots\ \Xi_{J}.
\ee
Here the size of $\Xi_{1}$ has been cut down a bit to $\tilde\Xi_{1}$ to make room for the set $\aleph$. 

The ``correct'' definition of $S_{N}(\gt)$ is then to replace \eqref{eq:SNbad} by:
\be\label{eq:SNdef}
S_{N}(\gt):=\sum_{\g\in\gW_{N}}e(\gt\<\g e_{2},e_{2}\>)
,
\ee
from which the argument follows as described above. In the end, we prove Theorem \ref{thm:mainN}, and hence Theorem \ref{thm:main}.
As already mentioned, Theorem \ref{thm:primes} is an easy corollary to Theorem \ref{thm:mainN}.

\begin{rmk}\label{rmk:McMullen}
One may  ponder the flexibility of our methods in 
applications
to other problems. For one in particular, McMullen \cite{McMullen2009} has
popularized the problem of producing many closed geodesics in a compact subset of the modular surface, defined over a fixed real quadratic number field $\Q(\sqrt{f})$. This is the same as producing many elements  $\g\in\cG$ so that $\tr(\g)^{2}-4$ has square-free part $f$
. 
Specifically, McMullen asks whether there is a finite alphabet $\cA$ so that the set of traces in $\cG_{\cA}$ contains every sufficiently large admissible integer.
Our use of 
Vinogradov's bilinear estimates relies crucially on the structure in \eqref{eq:Ge2e2} and does not apply as it stands to the problem of traces. 
We plan to return to this problem in the future.
\end{rmk}

The proof of Theorem \ref{thm:gdPlus} follows along completely different lines, and is 
inspired by the recent advances in projection theorems \cite{Bourgain2010}. The observation here is that the set $\fD$ has a certain ``sum-set'' structure. Namely, if $b/d\in\fR$ is a reduced fraction and $a\in\cA$, then clearly
\be\label{eq:sumSet}
\frac1{a+\frac bd}=\frac d{b+ad}\in\fR.
\ee
This implies that $b+ad\in\fD$ whenever $b/d\in\fR$ and $a\in\cA$; we exploit this sum-set structure to produce the bound \eqref{eq:gdPlus}.
We note further that the Discretized Ring Theorem \cite{Bourgain2003} can be used to get an exponent gain over the lower bound \eqref{eq:DANpower} even when $\gd\le 1/2$. 


\subsection{Outline of the Paper}

In \S\ref{sec:expVect}, we study the multiplicative properties of expanding  eigen-values and -vectors for matrices in $\G$.
We use \S\ref{sec:Setup}
to construct the main ensemble $\gW_{N}$, reserving the 
construction of the leading set $\aleph$ for  \S\ref{sec:New}.
The major arc analysis is carried out in \S\ref{sec:Maj}, while the minor arc bounds are proved in \S\S\ref{sec:Bnd1}--\ref{sec:Bnd2}.  Theorem \ref{thm:mainN} is then proved in \S\ref{sec:posProp}, as is its corollary, Theorem \ref{thm:primes}. Lastly, we prove Theorem \ref{thm:gdPlus} in \S\ref{sec:appendix}.


\subsection*{Notation}

Throughout we use the following standard notation. We write $f\sim g$ to mean $f/g\to1$. We use the Landau/Vinogradov notations $f=O(g)$ and $f\ll g$ synonymously  to mean there exists an implied constant $C>0$ such that for $x$ sufficiently large, $f(x)\le C g(x)$. Moreover $f\asymp g$ denotes $f\ll g\ll f$.
We allow the implied constants to depend at most  on the fixed alphabet $\cA$, unless otherwise specified. We also use the  short hand $e(x)=e^{2\pi ix}$. 
The cardinality of a finite set $S$ is denoted both as $\#S$ and $|S|$, and the Lebesgue measure of an interval $\cI$ is also $|\cI|$.
Throughout there are some constants $c,C>0$ which may change from line to line. 

\subsection*{Acknowledgements} We thank Curt McMullen for bringing this problem to our attention, and Doug Hensley 
and Peter Sarnak
for many helpful comments and suggestions regarding this work.



\section{
Large Matrix Products}
\label{sec:expVect}


In this section, we develop some tools in large matrix products, reminiscent of the  avalanche principle, see e.g. \cite[Ch. 6]{Bourgain2005} or \cite[\S2]{GoldsteinSchlag2001}.
%
Recall that
 $\G=\G_{\cA}$ is 
  the semigroup generated by even words in the matrices
\eqref{eq:cGgens}, for
$a\in \cA$.
An easy induction shows that for
$
\g=\bigl( \begin{smallmatrix}
a&b\\ c&d
\end{smallmatrix} \bigr)
\in\G,
$
$\g\neq I$,
we have 
$$
1\le a\le \min(b,c)\le\max(b,c)< d.
$$

We use the Frobenius norm:
\be\label{eq:supToFrob}
\|\g\|
:=\sqrt{a^{2}+b^{2}+c^{2}+d^{2}}
.
\ee
Note 
 that the trace and norm are comparable up to constants:
\be\label{eq:normToTrace}
\foh\|\g\|
\le \tr\g\le 
2\|\g\|
,
\ee
as are the norm, sup-norm, and ``second column'' norm:
\be\label{eq:normToE2}
\|\g\|_{\infty}=d<
|\g e_{2}|
=
\sqrt{b^{2}+d^{2}}
<
\|\g\|
<
2|\g e_{2}|
<
4\|\g\|_{\infty}
.
\ee

For $\g\in\G$, let the expanding and contracting eigenvalues of $\g$ be $\gl_{+}(\g)$ and $\gl_{-}(\g)=1/\gl_{+}(\g)$, with corresponding normalized eigenvectors $v_{+}(\g)$ and $v_{-}(\g)$.
Write $\gl=\gl_{+}$ for the expanding eigenvalue, so that
$$
\gl(\g)
=\gl_{+}(\g)
=
{\tr(\g)+ \sqrt{\tr(\g)^2 - 4 
}\over2}
.
$$
Note that for all $\g\in\G$,
the eigenvalues are real, and $\gl>1$ for $\g\neq I$. 
We require the following elementary but very useful observation.


\begin{prop}
The
eigenvalues of 
two 
matrices
$\g,\g'\in\G$ with large norms 
behave essentially multiplicatively,
subject to their  expanding eigenvectors   facing nearby directions. 
That is,
\be\label{eq:1.7}
\gl(\g\g')
=
 \gl(\g)\gl(\g')
\left[
1
+O\left(
\big|v_{+}(\g)-v_{+}(\g')\big|
+\frac1{\|\g\|^{2}}+\frac1{\|\g'\|^{2}}
\right)
\right]
.
\ee
%

Moreover, the expanding eigenvector of the product $\g\g'$ faces 
a nearby 
direction to that of the first $\g$, (and the same in reverse),
\be\label{eq:1.8}
|v_{+}(\g\g') - v_{+}(\g)|\ll \frac1{\|\g\|^{2}}
\qquad
\text{ and }
\qquad
|v_{-}(\g\g') - v_{-}(\g')|\ll \frac1{\|\g'\|^{2}}
.
\ee
The  implied constants above are absolute.
\end{prop}

\pf


For $\g$ large, we have:
\be
\label{eq:glToTr}
\gl(\g)
=
{\tr(\g)+ \sqrt{\tr(\g)^2 - 4 
}\over2}
=
\tr(\g)+
O\left(\frac1{\|\g\|}\right)
,
\ee
and

\beann
v_{+}(\g)
&=&
{
(b,\gl_{+}(\g)-a)
\over
\sqrt{b^{2}+(\gl_{+}(\g)-a)^{2}}
}
=
{
(b,d)
\over
\sqrt{b^{2}+d^{2}}
}
+O\left(\frac1{\|\g\|^{2}}\right)
,
\\
\nonumber
v_{-}(\g)
&=&
{
(d-\gl_{-}(\g),c)
\over
\sqrt{(d-\gl_{-}(\g))^{2}+c^{2}}
}
=
{
(-d,c)
\over
\sqrt{c^{2}+d^{2}}
}
+O\left(\frac1{\|\g\|^{2}}\right)
.
\eeann
Note that for $\g$ large,
\be\label{eq:vpvmperp}
|\<v_{+}(\g), v_{-}(\g)^{\perp}\>| = {bc+d^{2}\over \sqrt{b^{2}+d^{2}}\sqrt{c^{2}+d^{2}}} + O\left(\frac1{\|\g\|^{2}}\right) \ge \frac12,
\ee
meaning that the angle between expanding and contracting vectors does not degenerate. 

By \eqref{eq:glToTr}, it is enough to show that the traces behave essentially multiplicatively.
We compute:
\beann
|\tr(\g\g')-\tr(\g)\tr(\g')|
&=&
|(aa'+bc'+cb'+dd') - (a+d)(a'+d')|
\\
&\le&
\frac d{d'}
\left|
{bc'd'\over d}-a'd'
\right|
+
\frac{d'}d
\left|{cb'd \over d'}-ad\right|
\\
&\le&
\frac d{d'}
\left(
1 +
c'\left|
{bd'\over d}-b'
\right|
\right)
+
\frac{d'}d
\left(
1 +
c\left|
{b'd\over d'}-b
\right|
\right)
\\
&=&
\frac d{d'}
+
\frac{d'}d
+
(
cd'+
c' d
)
\left|
{b\over d}-{b'\over d'}
\right|
.
\eeann
We clearly have
\beann
\left|
\frac bd
-
\frac{b'}{d'}
\right|
&=&
|v_{+}(\g)-v_{+}(\g')|
+
O\left(\frac1{\|\g\|^{2}}
+
\frac1{\|\g'\|^{2}}\right)
,
\eeann
and hence
\beann
|\tr(\g\g')-\tr(\g)\tr(\g')|
&\ll&
d
d'
\left(
|v_{+}(\g)-v_{+}(\g')|
+
\frac1{\|\g\|^{2}}
+
\frac1{\|\g'\|^{2}}
\right)
.
\eeann
From this and \eqref{eq:glToTr}, \eqref{eq:1.7} follows easily.
One proves \eqref{eq:1.8} in a similar fashion.
%
\epf


\section{Construction of $\gW_{N}$}\label{sec:Setup}

\subsection{The leading term $\aleph$}\label{sec:set0}\

In this subsection, we posit the existence and all necessary properties of the leading set $\aleph$ 
used
in our construction of the main ensemble $\gW_{N}$. The proof of its existence is 
arguably the most
 technical part of the whole paper, so in the interest of exposition, we postpone it to \S\ref{sec:New}. 

Once and for all, we fix a density point $x\in\fC$, and let
\be\label{eq:vDef}
\fv=\fv_x:={(x,1)\over\sqrt{1+x^{2}}}
\ee
be the corresponding unit vector. We will henceforth be largely concerned with elements 
 of $\G$ whose expanding eigenvectors point in this direction.

For ease of exposition, 
we assume henceforth that for all $q\ge1$, the reduction of $\G$ is full,
\be\label{eq:GamToSL2q}
\G(\mod q)\cong\SL_{2}(q),
\ee 
which is anyway the case for any alphabet $\cA$ containing $1$ and $2$; see Remark \ref{rmk:StongApp}.
Minor modifications are needed in the general case.

For $N$ large 
and $\gd$ 
exceeding $\gd_{0}$ in \eqref{eq:gdAllow}%
, let
\be\label{eq:gbIs}
\Bpow:=
\frac1{1000}
(
\gd-\gd_{0}
)
>0
,
\ee
and let $\ga_{0}>0$ be a parameter to be chosen later in  \eqref{eq:ga0Is}.
Then we set
\be\label{eq:Bis}
B
:=
N^{\Bpow}
,
\ee
and 
\be\label{eq:cQIs}
\cQ:=N^{
\ga_{0}
/ \log\log N}
.
\ee 

Let 
\be\label{eq:cUrange}
\cU
\subset\left[\cIm B,\cIp B\right]
\ee 
be an arithmetic progression of real numbers starting with $u_{0}=\cIm B$ having common difference 
\be\label{eq:ujDiff}
|u-u'|=2B/\cQ^{5}
,
\ee
for $u,u'$ consecutive terms in $\cU$, and ending with $u>(\cIp -{2\over\cQ^{5}})B$.
Then the cardinality of $\cU$ is
\be\label{eq:cUsize}
|\cU|\asymp \cQ^{5}
.
\ee

\begin{prop}\label{prop:cBII}
For each $u\in\cU$, there are non-empty  sets
$\aleph_{u}\subset\G$, 
all of the same cardinality
\be\label{eq:alephUsize}
|\aleph_{u}|=|\aleph_{u'}|
,
\ee
so that the following holds.
For every $\fa\in\aleph_{u}$,
its expanding eigenvector is restricted by
\be\label{eq:d8p}
|v_{+}(\fa)-\fv|<\cQ^{-5}
,
\ee
and
 its expanding eigenvalue $\gl(\fa)$ is restricted by
\be\label{eq:d9p}
|\gl(\fa)-u|<{B\over \cQ^{5}}
.
\ee
In particular,
\be\label{eq:glfaToB}
\frac1{200}B<\gl(\fa)<B,
\ee 
for $N$ large.
Moreover, for any
$q<\cQ$, any
$\gw\in\SL_{2}(q)$ and any
$u\in\cU$, we have
\be\label{eq:cBSatp}
\#
\left\{
\fa\in\aleph_{u}
:
\fa\equiv \gw(\mod q)
\right\}
=
 {|\aleph_{u}|
\over  |\SL_{2}(q)|}
(
1+
O(
\cQ^{-4})
)
,
\ee
where the implied constant does not depend on $q$, $\gw$, or $u$.
\end{prop}

With the sets $\aleph_{u}$ as above, we define the main leading set $\aleph$ to 
be the union of the sets $\aleph_{u}$,
\be\label{eq:alephDecomp}
\aleph:=\bigsqcup_{u\in\cU}\aleph_{u}
\ee
Note that the sets $\aleph_{u}$ are disjoint by \eqref{eq:d9p} and \eqref{eq:ujDiff}.
We repeat that the proof of Proposition \ref{prop:cBII} will be postponed to \S\ref{sec:New}.

\subsection{Sector Counting}\label{sec:Lalley}\

In this section we give the following slight refinement of Hensley's estimate \eqref{eq:2dim}, which follows directly from Lalley's methods \cite{Lalley1989}.
\begin{prop}\label{prop:Lalley}
There is a constant $\fc=\fc(\cA)>0$ so that as long as $H<T^{\fc/\log\log T}$, we have
\be\label{eq:Lalley}
\#\left\{\g\in\G:\|\g\|<T\text{ and }|v_{+}(\g)-\fv|<\frac1H\right\}
\gg
{T^{2\gd}
\over H}
,
\ee
as $T\to\infty$. 
\end{prop}
\pf[Sketch of proof]
Lalley \cite[Theorem 9]{Lalley1989}
proves the asymptotic formula 
\be\label{eq:Lalley1}
\text{Left-hand side of }\eqref{eq:Lalley}\sim C\cdot T^{2\gd}\mu(\cI)
\ee
under the assumption that $\G$ is a non-elementary convex-cocompact sub{\it group} of $\SL_{2}(\R)$. 
Here $\cI$ is the interval of length $1/H$ about $\fv$, and $\mu$ is the $\gd$-dimensional Hausdorff measure 
supported on the limit set $\fC$, lifted (by abuse of notation) to $\bP^{1}$ via 
$$
d\mu(x,y)=d\mu(x/y),
$$
$y\neq0$.
After setting up the symbolic dynamics, the requirement that $\G$ not contain parabolic elements is needed in the renewal method to make the distortion function eventually positive, see \cite[pp. 33, 41] {Lalley1989}. Our semigroup $\G$ has no parabolic elements, so the only difference here between a group and (free) semigroup is that for the latter, the transition matrix (see \cite[pp. 5, 32]{Lalley1989}) is trivial, that is, all sequences are allowed in the symbolic dynamics.
The rate in \eqref{eq:Lalley1} can be determined directly from Lalley's method (see \cite[\S12]{BourgainGamburdSarnak2011}), with the error crudely estimated as 
\be\label{eq:LalleyError}
\ll T^{2\gd-c/\log\log T}.
\ee 
Since $\fv$ in \eqref{eq:vDef} corresponds to a density point in $\fC$, we have, again crudely, that 
$$
\mu(\cI)\gg_{\vep} H^{-\gd-\vep}\gg H^{-1},
$$ 
since $\gd<1$. A sufficient condition for the main term, being bounded below by $CT^{2\gd}/H$, to dominate the error in \eqref{eq:LalleyError}, is that $H<T^{\fc/\log\log T}$ with $\fc<c$.
\epf

\begin{rmk}
The methods of Dolgopyat \cite{Dolgopyat1998} and Naud \cite{Naud2005} could be used to prove \eqref{eq:Lalley} with $H$ as large as $T^{\vep}$, but this is not needed in our applications. 
\end{rmk}

With this crude estimate in hand, we proceed in the next subsection to detail our  construction of the special sets $\Xi$ alluded to in \S\ref{sec:outline}.

\subsection{The set $\Xi(M,H;L,k)$}\label{sec:set1}\


\begin{prop}\label{prop:Xi}
Given $M\gg1$ and $H<M^{\fc/\log\log M}$, there exists some $L$ in the range
\be\label{eq:LtoMbnd}
\frac14M \le L \le 4M
,
\ee
an integer $k\asymp\log M$, and a set $\Xi=\Xi(M,H;L,k)\subset\G$ having the following properties.
For all $\g\in\Xi$, 
the expanding eigenvalues are controlled to within $1/\log L$:
\be\label{eq:4.0}
L\left(1-\frac1{\log L}\right)<\gl(\g)<L,
\ee
the expanding eigenvectors are controlled to within $1/H$:
\be\label{eq:4.0'}
|v_{+}(\g)-\fv|<\frac1{H},
\ee
and the wordlength  metric $\ell$  (in the generators \eqref{eq:Ggens} of $\G$) is controlled exactly:
\be\label{eq:wdlngth}
\ell(\g)=k.
\ee
Moreover, the cardinality of $\Xi$ is controlled by
\be\label{eq:XiSize}
L^{2\gd}\gg \#\Xi\gg {L^{2\gd}\over H (\log L)^{2}}
.
\ee
\end{prop}

Recall again the the implied constants depend at most on $\cA$, which is thought of as 
fixed throughout.

\pf

We proceed by the following algorithm.

\begin{enumerate}

\item
Let  $S_{1}\subset \G$ be the set  of $\g\in\G$ of norm controlled by
$
\|\g\|\asymp M
$
and
for which the expanding vector $v_{+}(\g)$ is within $\frac1{H}$ of the fixed vector  $\fv$:
$$
S_{1}:=\left\{\g\in\G:\frac M2<\|\g\|<M,|v_{+}(\g)-\fv|<\frac1H\right\}.
$$
By \eqref{eq:Lalley}, we have that 
 $$
 \#S_{1}\gg {M^{2\gd}\over H}.
 $$

 \item
By \eqref{eq:normToTrace} and \eqref{eq:glToTr},
expanding eigenvalues $\gl(\g)$ of $\g\in S_{1}$ satisfy
$$
\frac14M \le \gl(\g) \le 4M
.
$$
Hence we can find (by pigeonhole) an $L$ in this range so that 
$$
\#\{\g\in S_{1}: L\left(1-\frac1{\log L}\right)<\gl(\g)< L\}
\gg
{L^{2\gd}\over H\log L}
.
$$
Call the above set $S_{2}$; its expanding eigenvalues are all nearly of the same size.
\item
Lastly, note that the wordlength metric $\ell$ is commensurable with the archimedean one, 
$$
\ell(\g)\asymp \log \|\g\|,
$$
with implied constant depending on $\cA$. So (again by pigeonhole) we can find some $k$ such that
\be\label{eq:wdlength}
\#\{\g\in S_{2}: \ell(\g)=k
\}
\gg
{L^{2\gd}\over H(\log L)^{2}}
.
\ee
Call this set $S_{3}$;
then the elements of $S_{3}$ all have the same wordlength, in addition to the previous qualities.
\end{enumerate}

We rename this last set $S_{3}$ to $\Xi
=\Xi(M,H;L,k)$.
\epf

\subsection{Decomposing $N$ and the ensemble $\gW_{N}$}\label{sec:gWN}\

We return to our main parameter $N$,
and decompose it 
dyadically
as follows.
Recall that we have already presupposed the construction of a set $\aleph$ in \eqref{eq:alephDecomp}, all of whose expanding eigenvectors are within $\cQ^{-5}$ of $\fv$, and with eigenvalues of size $B$, see \eqref{eq:glfaToB}. Recall from  \eqref{eq:Bis} that $B=N^{\Bpow}$, and that $\cQ$ is given by \eqref{eq:cQIs}.


{\bf Setup:}
We start by taking 
\be\label{eq:66s}
M=\sqrt N/B=N^{1/2-\Bpow}, \qquad H=\cQ^{5}.
\ee
The exponent $\ga_{0}$ in the definition \eqref{eq:cQIs} of $\cQ$
 is chosen in \eqref{eq:ga0Is} to be
 sufficiently small that $H<M^{\fc/\log\log M}$. Run the algorithm of the previous subsection to generate the set $\Xi(M,H;L,k)$. By \eqref{eq:LtoMbnd}, the returned parameter $L$ satisfies 
$$
L=\ga_{1}M=\ga_{1}N^{1/2-\Bpow},
$$
with
$$
\ga_{1}\in(1/4,4).
$$
Write 
$$
\tilde N_{1}:=L=\ga_{1}N^{1/2-\Bpow},\qquad
N_{1}:=\ga_{1}N^{1/2}=B\cdot \tilde N_{1}
,
$$ 
and rename the returned set to $\tilde\Xi_{1}=\Xi(M,H;L,k)$, also setting
$$
\Xi_{1}:=\aleph\cdot\tilde\Xi_{1}
.
$$

\begin{rmk}\label{rmk:wdXi}
Despite the wordlength in $\aleph$ being unrestricted, the wordlength in $\tilde\Xi_{1}$ is fixed in \eqref{eq:wdlngth}. So the representation of an element in $\Xi_{1}$ as  a product of ones in $\aleph$ and $\tilde\Xi_{1}$ is still unique.
\end{rmk}

We have crudely that 
\be\label{eq:XiSize'}
|\Xi_{1}|\ge |\tilde\Xi_{1}|\gg_{\vep} \tilde N_{1}^{2\gd-\vep}\gg N^{\gd-2\gd\Bpow-\vep}.
\ee
(The cardinality of $\aleph$ is quite deficient relative to its norm, so we lose little from estimating trivially $|\aleph|\ge1$.)

{\bf Step 1:}
Next we set 
$$
M=
{N_{1}^{1/2}\over \ga_{1}}
=
{N^{1/4}\over \ga_{1}^{1/2}}
,
\qquad
H=\log M,
$$ 
and generate another set $\Xi(M,H;L,k)$. Define 
$$
N_{2}:=L = \ga_{2}M = {\ga_{2} N^{1/4}\over \ga_{1}^{1/2}}
,
$$
with $\ga_{2}\in(1/4,4)$, and 
rename the returned set to $\Xi_{2}.$ 
We have
$$
|\Xi_{2}|\gg {N_{2}^{2\gd}\over (\log N_{2})^{3}}
.
$$

{\bf Iterate:}
Start with $j=3$ and iterate up to $j= J-1$, where 
\be\label{eq:JcIs}
2^{J-1}=c\log N.
\ee 
Here the constant $c>0$ is absolute (independent of $N$), determined by \eqref{eq:NJis}. 
For each such $j$, set
\be\label{eq:66ss}
M:=
{
(N_{j-1})^{1/2}
\over 
\ga_{j-1}
}
=
{N^{1/2^{j}}\over \ga_{j-1}^{1/2} \ga_{j-2}^{1/4}\cdots (\ga_{1})^{1/2^{(j-1)}}}
,
\qquad
H=\log M,
\ee
and use Proposition \ref{prop:Xi} generate the set $\Xi(M,H;L,k)$. Define 
\be\label{eq:66s4}
N_{j}:=L=\ga_{j}M={\ga_{j}N^{1/2^{j}}\over \ga_{j-1}^{1/2} \ga_{j-2}^{1/4}\cdots (\ga_{1})^{1/2^{(j-1)}}},
\ee
with $\ga_{j}\in(1/4,4)$, and
call the returned set $\Xi_{j}.$ Note that
\be\label{eq:XijSize}
|\Xi_{j}|\gg {N_{j}^{2\gd}\over (\log N_{j})^{3}}
\ee
and
\be\label{eq:NjToN2j}
\frac1{16}N^{1/2^{j}}<N_{j}<16 N^{1/2^{j}}
.
\ee

{\bf  End:}
For the last step,  $j=J$, we set
$$
M=
{ N_{J-1}
\over 
(\ga_{J-1})^{2}}
=
  {N^{1/2^{(J-1)}}\over\ga_{J-1}(\ga_{J-2})^{1/2}\cdots\ga_{1}^{1/2^{(J-2)}}},
  \qquad
  H=\log M,
$$ 
and generate one last set $\Xi_{J}:=\Xi(M,H;L,k)$. Define 
$$
N_{J}:=L =  {\ga_{J}N^{1/2^{(J-1)}}\over\ga_{J-1}\cdots\ga_{1}^{1/2^{(J-2)}}}
\asymp
N^{1/2^{(J-1)}}
=
e^{1/c}
\ll1
,
$$ 
where we used \eqref{eq:JcIs}.
Since $\frac14<N_{J}/ M
=\ga_{J}<4$, we have
\be\label{eq:prodNJtoN}
\frac14
<
{N_{1}N_{2}\dots N_{J}\over N}
=
{B\tilde N_{1}N_{2}\dots N_{J}\over N}
<4
.
\ee




We now define the main ensemble $\gW_{N}$ by concatenating the sets $\Xi_{j}$ developed above. 
\be\label{eq:gWLdef}
\gW_{N}:
=
\Xi_{1}\cdot \Xi_{2}\cdots\Xi_{J-1}\cdot \Xi_{J}
=
\aleph\cdot \tilde\Xi_{1}\cdot \Xi_{2}\cdots\Xi_{J-1}\cdot \Xi_{J}
.
\ee

\subsection{Properties of  $\gW_{N}$}\

For $\g\in\gW_{N}$, write
$$
\g=\fa\cdot\tilde\xi_{1}\ \xi_{2}\cdots\xi_{J}
$$
according to the decomposition \eqref{eq:gWLdef}. Note that by the fixed wordlength restriction \eqref{eq:wdlngth}, this decomposition is unique, see Remark \ref{rmk:wdXi}.
Recall that the expanding vectors $v_{+}$ all point nearly in the direction of $\fv$ in \eqref{eq:vDef}.

\begin{lem}\label{lem:5.6}
For any
$2\le j_{1}\le j_{2}\le J$,
and 
$\xi_{j_{1}}\in\Xi_{j_{1}},\cdots \xi_{j_{2}}\in\Xi_{j_{2}}$,
and any $\fa\in\aleph$, $\tilde\xi_{1}\in\tilde\Xi_{1}$,
we have the 
following control on expanding eigen-vectors and -values of large products:
\be\label{eq:vpXi}
|v_{+}(\tilde\xi_{1}\cdot\xi_{2}\cdots\xi_{J})
-\fv|
\ll
\cQ^{-5}
,
\ee
\be\label{eq:66sss}
\frac12<{\gl(\xi_{j_{1}}\xi_{j_{1}+1}\cdots \xi_{j_{2}-1}\xi_{j_{2}})\over N_{j_{1}}N_{j_{1}+1}\cdots N_{j_{2}-1} N_{j_{2}}}<2
,
\ee
\be\label{eq:66s8}
\frac12<{\gl(\tilde \xi_{{1}}\xi_{2}\cdots \xi_{j_{2}-1}\xi_{j_{2}})\over \tilde N_{{1}}N_{2}\cdots N_{j_{2}-1} N_{j_{2}}}<2
,
\ee
and
\be\label{eq:66s9}
\frac12<{\gl(\fa\tilde \xi_{{1}}\xi_{2}\cdots \xi_{j_{2}-1}\xi_{j_{2}})\over \gl(\fa)\tilde N_{{1}}\cdot N_{2}\cdots N_{j_{2}-1} N_{j_{2}}}<2
.
\ee
\end{lem}

\pf
From \eqref{eq:1.8}, \eqref{eq:4.0'}, and the choice of $H$ in \eqref{eq:66s}, we have that
\beann
|v_{+}(\tilde \xi_{1}\cdot\xi_{2}\cdots\xi_{J})-\fv|
&\le&
|v_{+}(\tilde \xi_{1}\cdot\xi_{2}\cdots\xi_{J})-v_{+}(\tilde \xi_{1})|
+
|v_{+}(\tilde \xi_{1})-\fv|
\\
&\ll&
\frac1{\|\tilde\xi_{1}\|^{2}}
+
\frac1{\cQ^{5}}
,
\eeann
whence \eqref{eq:vpXi} follows from \eqref{eq:cQIs}.

Similarly, we have for $j\in[j_{1},j_{2}]\subset[2,J]$ that
\be\label{eq:66s5}
|v_{+}(\xi_{j}\xi_{j+1}\cdots \xi_{j_{2}})-\fv|
\ll
\frac1{\log N_{j}}
,
\ee
where we used the choice of $H$ in \eqref{eq:66ss}.

We now prove by downward induction on $j_{1}$ that
\bea
\label{eq:66s7}
\gl(\xi_{j_{1}}\xi_{j_{1}+1}\cdots\xi_{j_{2}})
&=&
N_{j_{1}}N_{j_{1}+1}\cdots N_{j_{2}}
\\
\nonumber
&&
\times
\left[
1+
O
\left(
\frac1{\log N_{j_{1}}}
+
\frac1{\log N_{j_{1}+1}}
+\cdots +
\frac1{\log N_{j_{2}}}
\right)
\right]
.
\eea
If $j_{1}=j_{2}$, then  \eqref{eq:66s7} follows immediately from \eqref{eq:4.0} and \eqref{eq:66s4}. If $j_{1}=j_{2}-1$, then from \eqref{eq:1.7}, \eqref{eq:4.0}, \eqref{eq:4.0'}, and \eqref{eq:66s4}, we have
\beann
\gl(\xi_{j_{2}-1}\xi_{j_{2}})
&=&
\gl(\xi_{j_{2}-1})\gl(\xi_{j_{2}})
\\
&&
\times
\left[
1+
O
\left(
|v_{+}(\xi_{j_{2}-1})-v_{+}(\xi_{j_{2}})|
+
\frac1{\|\xi_{j_{2}-1}\|^{2}}
+
\frac1{\|\xi_{j_{2}}\|^{2}}
\right)
\right]
\\
&=&
N_{j_{2}-1}N_{j_{2}}
\left[
1+
O
\left(
\frac1{\log N_{j_{2}-1}}
+
\frac1{\log N_{j_{2}}}
\right)
\right]
,
\eeann
as desired.

In general, we have by \eqref{eq:66s5} that
\beann
\gl(\xi_{j_{1}}\xi_{j_{1}+1}\cdots\xi_{j_{2}})
&=&	
\gl(\xi_{j_{1}})\gl(\xi_{j_{1}+1}\cdots\xi_{j_{2}})
\\
&&
\hskip-1in
\times
\left[
1+
O
\left(
|v_{+}(\xi_{j_{1}})-v_{+}(\xi_{j_{1}+1}\cdots\xi_{j_{2}})|
+
\frac1{\|\xi_{j_{1}}\|^{2}}
+
\frac1{\gl(\xi_{j_{1}+1}\cdots\xi_{j_{2}})^{2}}
\right)
\right]
\\
&=&	
N_{j_{1}}\gl(\xi_{j_{1}+1}\cdots\xi_{j_{2}})
\left[
1+
O
\left(
\frac1{\log N_{j_{1}}}
+
\frac1{\log N_{j_{1}+1}}
\right)
\right]
,
\eeann
from which \eqref{eq:66s7} follows by induction.

The rate in \eqref{eq:66s7} may be replaced crudely by 
\be\label{eq:NJis}
\left[
1+
O
\left(
{2^{J}\over \log N}
\right)
\right]
,
\ee
whence \eqref{eq:66sss} follows on taking the constant $c$ in \eqref{eq:JcIs} sufficiently small (independent of $N$). The estimates \eqref{eq:66s8} and  \eqref{eq:66s9} are proved in the same way.
\epf

As a consequence of \eqref{eq:66s9}, \eqref{eq:glfaToB}, and \eqref{eq:prodNJtoN},  we have that for all $\g\in\gW_{N}$,
\be\label{eq:normGsize}
\|\g\|\le 2 \gl(\g) \le 16N
,
\ee
so indeed the norms are all controlled.

Moreover the size of $\gW_{N}$ is not too much smaller than \eqref{eq:SN0}.
Indeed, we have from \eqref{eq:XiSize'}, \eqref{eq:XijSize}, and  \eqref{eq:JcIs} that 
\bea
\nonumber
\#\gW_{N}
&=&
\#\Xi_{1}\cdot
\#\Xi_{2}\cdot
\cdots
\#\Xi_{J}
\\
\nonumber
&\gg_{\vep}&
\tilde N_{1}^{2\gd-\vep}
{(N_{2})^{2\gd}\over (\log N_{2})^{3}}
\cdots
{(N_{J})^{2\gd}\over (\log N_{J})^{3}}
\\
\label{eq:gWLsize}
&\gg&
{N^{2\gd-2\gd\Bpow-\vep}
}
.
\eea
%
It also follows that
for any $j\ge2$, 
\be\label{eq:GWLsize2}
\#\Xi_{j}\cdot
\#\Xi_{j+1}\cdot
\cdots
\#\Xi_{J}
\gg
(N_{j}N_{j+1}\cdots N_{J})^{2\gd}
\cdot
e^{-c(J-j)\log\log N_{j}}
,
\ee
for an absolute constant $c>0$.

With the set $\gW_{N}$ constructed, we define our exponential sum $S_{N}$ as in \eqref{eq:SNdef}, and proceed with the circle method.




\section{Major Arcs Analysis}\label{sec:Maj}

 In this section we
estimate the
  major arcs contribution.
First we use the set $\aleph$ 
described
in
\S\ref{sec:set0}  to prove that in the major arcs, 
our exponential sum
$S_{N}$ 
in \eqref{eq:SNdef}
splits as a product of modular and archimedean components, as in \eqref{eq:SNvqvpgb}. Then we prove that the major arcs contribution is of the correct order of magnitude.

\subsection{Splitting into Modular and Archimedean Components}\

Let 
$\cQ
$ 
%
be as in \eqref{eq:cQIs}, and $B$ as in \eqref{eq:Bis}.
Recall from \eqref{eq:MajArcsDef} that the major arcs of level $\cQ$ are given by
$$
\fM_{\cQ}=
\bigsqcup_{q<\cQ}\bigsqcup_{(a,q)=1}\left[\frac aq-\frac{\cQ}N,\frac aq+\frac{\cQ}N\right].
$$
Let $\nu_{q}:\Z/q\Z\to\C$ 
record the
mod $q$
 distribution of $\fD$. 
That is, for $a\in\Z/q\Z$, set
\be\label{eq:nuQdef}
\nu_{q}(a)
:=
\frac1{|\SL_{2}(q)|}
\sum_{\gw\in\SL_{2}(q)}
e\left(\frac aq\<\gw e_{2},e_{2}\>\right)
.
\ee

\begin{thm}\label{thm:split}
There exists a function
$
\vp_{N}:\R/\Z\to\C
,
$
given explicitly in \eqref{eq:varpiIs},
 satisfying 
the following three conditions. 
\begin{enumerate}
\item\label{it:1} The Fourier transform 
$$
\hat\vp_{N}:\Z\to\C:n\mapsto\int_{0}^{1}\vp_{N}(\gt)e(-n\gt)d\gt
$$ 
is 
real-valued and non-negative, with 
 \be\label{eq:vpSize}
\vp_{N}(0)=
\sum_{n} \hat\vp_{N}(n)\ll|\gW_{N}|
.
 \ee

\item\label{it:2} 
For
$\gbaa
 N< n< \gbbb
 N
,
$
we have
\be\label{eq:vpBndBelow}
\hat\vp_{N}(n)\gg {|\gW_{N}|\over N}
.
\ee

\item\label{it:3} 
Moreover, 
 we have
on the major arcs $\gt=\frac aq+\gb\in\fM_{\cQ}$ that
\be\label{eq:SNsplit}
S_{N}\left(\frac aq+\gb\right)
=
\nu_{q}(a)
\vp_{N}(\gb)
\left(1
+
O(
\cQ^{-4}
)
\right)
.
\ee
\end{enumerate}
\end{thm}

\pf
We use the decomposition \eqref{eq:gWLdef} in the form
\be\label{eq:gWNdecomp2}
\gW_{N}=\aleph\cdot\gW'
\ee
with 
$$
\gW'=\tilde\Xi_{1}\Xi_{2}\cdots\Xi_{J}
,
$$
so that
$$
S_{N}(\gt)=\sum_{\fa\in\aleph}\sum_{\g\in\gW'}e(\gt\<\fa\g e_{2},e_{2}\>).
$$

For
$\fa\in\aleph$,
recall 
from
\eqref{eq:glfaToB}
that
we have
$
\gl(\fa)\asymp B,
$ 
and 
from \eqref{eq:d8p} and \eqref{eq:vpXi} that
\be\label{eq:vPlusTofvS}
|v_{+}(\fa)-\fv|<
\cQ^{-5}
,
\qquad
|v_{+}(\g)-\fv|\ll \cQ^{-5}
.
\ee
%

To use the key property \eqref{eq:d9p} 
in the construction of $\aleph$, 
we need to convert
 the expression $\<\fa\g e_{2},e_{2}\>$ 
 in $S_{N}$
into one involving $\gl(\fa)$. 

 We will make regular use of the following elementary formula: For any two linearly independent vectors $v_{+},v_{-}\in\R^{2}$, we can write any $w\in\R^{2}$ as
\be\label{eq:vpmPerps}
w=
{\<w,v_{-}^{\perp}\>\over \<v_{+},v_{-}^{\perp}\>}v_{+}
+
{\<w,v_{+}^{\perp}\>\over \<v_{-},v_{+}^{\perp}\>}v_{-}
.
\ee
Here $(x,y)^{\perp}=(-y,x)$.
Recalling \eqref{eq:vpvmperp}, it easily follows  that for a unit vector $w$ and 
any large
$\xi\in\G$,
\be\label{eq:gwTogl}
\xi w=
\gl(\xi)
{
\<w,v_{-}^{\perp}(\xi)\>
\over
\<v_{+}(\xi),v_{-}^{\perp}(\xi)\>
}
v_{+}(\xi)
\left(
1+
O\left(
\frac1{\|\xi\|^{2}}
\right)
\right)
,
\ee
whence
\be\label{eq:xie2e2}
\<\xi e_{2},e_{2}\>=
\gl(\xi)
{
\<e_{2},v_{-}^{\perp}(\xi)\>
\over
\<v_{+}(\xi),v_{-}^{\perp}(\xi)\>
}
\<v_{+}(\xi),e_{2}\>
\left(
1+
O\left(
\frac1{\|\xi\|^{2}}
\right)
\right)
.
\ee

Applied to our present situation, we have by \eqref{eq:vPlusTofvS} that
\be\label{eq:ge2e2s}
\<\g e_{2},e_{2}\>=
\gl(\g)
{
\<e_{2},v_{-}^{\perp}(\g)\>
\over
\<\fv,v_{-}^{\perp}(\g)\>
}
\<\fv,e_{2}\>
\left(
1+
O\left(
\frac1{\cQ^{5}}
\right)
\right)
,
\ee
and
\bea
\nonumber
\<\fa\g e_{2},e_{2}\>
&=&
\gl(\fa\g)
{
\<e_{2},v_{-}^{\perp}(\fa\g)\>
\over
\<v_{+}(\fa\g),v_{-}^{\perp}(\fa\g)\>
}
\<v_{+}(\fa\g),e_{2}\>
\left(
1+
O\left(
\frac1{N^{2}}
\right)
\right)
\\
\label{eq:ge2e2sp}
&=&
\gl(\fa)\gl(\g)
{
\<e_{2},v_{-}^{\perp}(\g)\>
\over
\<\fv,v_{-}^{\perp}(\g)\>
}
\<\fv,e_{2}\>
\left(
1+
O\left(
\frac1{\cQ^{5}}
\right)
\right)
,
\eea
where we also used  \eqref{eq:1.7} and \eqref{eq:1.8}.
%

%
Comparing \eqref{eq:ge2e2sp} and \eqref{eq:ge2e2s}, we have that
\be\label{eq:fbTogl}
\<\fa\g e_{2},e_{2}\>
=
\gl(\fa)
\<
\g e_{2}
,e_{2}\>
+
O
\left(
N/\cQ^{5}
\right)
.
\ee

For $\gt=\frac aq+\gb\in\fM_{\cQ}$ with $|\gb|<\cQ/N$, we insert \eqref{eq:fbTogl} into $S_{N}$, giving
\bea
\nonumber
S_{N}\left(\frac aq+\gb\right)
&=&
\sum_{\fa\in\aleph}
\sum_{\g\in\gW'}
e
\left(
\frac aq
\<\fa\g e_{2},e_{2}\>
\right)
e
\left(
\gb
\<\fa\g e_{2},e_{2}\>
\right)
\\
\nonumber
&=&
\sum_{\fa\in\aleph}
\sum_{\g\in\gW'}
e
\left(
\frac aq
\<\fa\g e_{2},e_{2}\>
\right)
e
\bigg(
\gb
\gl(\fa)\<\g e_{2},e_{2}\>
\bigg)
+O\bigg(\cQ^{-4}|\gW|\bigg)
\\
%
\nonumber
&=&
\sum_{\g\in\gW'}
\sum_{\gw\in\SL_{2}(q)}
e
\left(
\frac aq
\<\gw\g e_{2},e_{2}\>
\right)
\sum_{\fa\in\aleph\atop\fa\equiv\gw(\mod q)}
e
\left(
\gb
\gl(\fa)
\<\g e_{2},e_{2}\>
\right)
\\
\label{eq:SNspl}
&&
+O\bigg(\cQ^{-4}|\gW_{N}|\bigg)
,
\eea
where we decomposed the $\fa$ sum into  residue classes $\gw$ in $\SL_{2}(q)$.

Next using \eqref{eq:alephDecomp}, write the innermost sum above as
\bea
\label{eq:cUdecomp}
&&
\hskip-.5in
\sum_{\fa\in\aleph\atop\fa\equiv\gw(\mod q)}
e
\left(
\gl(\fa)
\gb
\<\g e_{2},e_{2}\>
\right)
\\
\nonumber
&=&
\sum_{u\in\cU}
\sum_{\fa\in\aleph_{u}\atop\fa\equiv\gw(\mod q)}
e
\left(
\gl(\fa)
\gb
\<\g e_{2},e_{2}\>
\right)
\\
\nonumber
&=&
\sum_{u\in\cU}
e
\left(
\gb
u
\<\g e_{2},e_{2}\>
\right)
\left(
\sum_{\fa\in\aleph_{u}\atop\fa\equiv\gw(\mod q)}
1
\right)
\left(
1+O(\cQ^{-4})
\right)
,
\eea
where we applied \eqref{eq:d9p}.

By \eqref{eq:cBSatp} and \eqref{eq:alephUsize} (that the cardinality of $\aleph_{u}$ is the same for all $u$), the innermost sum is 
\be\label{eq:fbEqui}
\sum_{\fa\in\aleph_{u}\atop\fa\equiv\gw(\mod q)}
1
=
{|\aleph|\over |\cU|\cdot|\SL_{2}(q)|}
\left(
1+
O(\cQ^{-4})
\right)
,
\ee
where 
the implied constant does not depend on $u$, $\gw$, or $q$.

Returning to \eqref{eq:SNspl},
inputting
\eqref{eq:fbEqui}
and 
\eqref{eq:cUdecomp}
gives
\bea
\nonumber
S_{N}\left(\frac aq+\gb\right)
&=&
{1\over |\SL_{2}(q)|}
\sum_{\gw\in\SL_{2}(q)}
e
\left(
\frac aq
\<\gw e_{2},e_{2}\>
\right)
{|\aleph|\over |\cU|}
\sum_{\g\in\gW'}
\sum_{u\in\cU}
e
\left(
\gb
u
\<\g e_{2},e_{2}\>
\right)
\\
\label{eq:SNspl2}
&&
\hskip2in
\times
\left(
1
+
O\left(
\cQ^{-4}
\right)
\right)
,
\eea
where we used the fact that the $\gw$ sum runs over all of $\SL_{2}(q)$, so is independent of $\g$.
Note that $S_{N}$ has already split into modular and archimedean components, with the first 
piece
being $\nu_{q}(a)$ as in \eqref{eq:nuQdef}.

We continue to massage the archimedean component.
Fix $\g$ and $u$. For any $m\in\Z$ with 
$$
|m-u\<\g e_{2},e_{2}\>|
\le
B\<\g e_{2},e_{2}\>/\cQ^{5},
$$ 
we clearly have
$$
e
\left(
\gb
u
\<\g e_{2},e_{2}\>
\right)
=
e
\left(
\gb
m
\right)
\left(
1
+O(\cQ^{-4})
\right)
,
$$
and
there are $2B\<\g e_{2},e_{2}\>/\cQ^{5} + O(1)$ integers $m$ in this range.
Hence
\be\label{eq:uToMs}
e
\left(
\gb
u
\<\g e_{2},e_{2}\>
\right)
=
{\cQ^{5}\over2B\<\g e_{2},e_{2}\>}
\sum_{m\in\Z\atop
\left|{m\over\<\g e_{2},e_{2}\>}-u\right|
\le
{B\over \cQ^{5}}
}
e
\left(
\gb
m
\right)
\left(
1+O(\cQ^{-4})
\right)
.
\ee

Reversing the $u$ and $m$ sums and inserting \eqref{eq:uToMs} into \eqref{eq:SNspl2}
gives
\beann
S_{N}\left(\frac aq+\gb\right)
&=&
\nu_{q}(a)
\vp_{N}(\gb)
\left(
1+
O\left(
\cQ^{-4}
\right)
\right)
,
\eeann
where
\be\label{eq:varpiIs}
\vp_{N}(\gb)
:=
{|\aleph|\over |\cU|}
\sum_{\g\in\gW'}
{\cQ^{5}\over2B\<\g e_{2},e_{2}\>}
\sum_{m\in\Z}
e
\left(
\gb
m
\right)
\sum_{u\in\cU}
\bo_{
\left|{m\over\<\g e_{2},e_{2}\>}-u\right|
\le
{B\over \cQ^{5}}
}
.
\ee

Hence \eqref{eq:SNsplit} is satisfied.
Also the Fourier transform
\be\label{eq:vpNhat}
\hat\vp_{N}(n)
=
{|\aleph|\over |\cU|}
\sum_{\g\in\gW'}
{\cQ^{5}\over2B\<\g e_{2},e_{2}\>}
\sum_{u\in\cU}
\bo_{
\left|{n\over\<\g e_{2},e_{2}\>}-u\right|
\le
{B\over \cQ^{5}}
}
\ee
is clearly real and non-negative, so \eqref{it:1} is satisfied.

Now, combining
\eqref{eq:normToTrace}, 
\eqref{eq:normToE2},
\eqref{eq:glToTr},
\eqref{eq:66s8}, and \eqref{eq:prodNJtoN}, we have
\be\label{eq:ge2e2Size}
\frac1{4}\frac NB<\<\g e_{2},e_{2}\><4 \frac NB,
\ee
and hence
for 
$
\gbaa
 N
< n
 < \gbbb
 N
,
$
we have, crudely, that
$$
\cIm B
< 
{n\over\<\g e_{2},e_{2}\>}
<
\cIp
B
.
$$
Hence
by 
the spacing in
\eqref{eq:ujDiff}
of 
 $u\in\cU$ in this range,
 the innermost sum in \eqref{eq:vpNhat}
 is  guaranteed to have at least one contribution, giving
$$
\hat\vp_{N}(n)
\gg
{|\aleph|\over |\cU|}
\sum_{\g\in\gW'}
{\cQ^{5}\over2B\<\g e_{2},e_{2}\>}
\gg
{
|\aleph||\gW'|
\over N}
=
{|\gW_{N}|\over N}
,
$$
where we used \eqref{eq:ge2e2Size}, \eqref{eq:cUsize}, and \eqref{eq:gWNdecomp2}.
So \eqref{eq:vpBndBelow} is satisfied, and 
the proof of Theorem \ref{thm:split} is complete.
\epf

\subsection{The Major Arcs Contribution}\

Equipped with \eqref{eq:SNsplit}, it is now straightforward to 
produce 
the necessary
major arcs contribution. 
For technical reasons, we need a smoothed cutoff, and introduce
the triangle function,
 $\psi
,$ given by
\be\label{eq:psiIs}
\psi(x):=\threecase{1+x,}{if $-1<x<0$,}{1-x,}{if $0\le x<1$,}{0,}{otherwise.}
\ee
It is well known 
 that 
the Fourier transform is non-negative:
\be\label{eq:psiHatIs}
\hat\psi(y)=\left({\sin(\pi y)\over \pi y}\right)^{2}.
\ee
Let $\psi_{N}$ be the function localized at level $\cQ/N$ near the origin:
$$
\psi_{N}(x)
:=
\psi\left({N\over \cQ}x\right)
$$
Periodize $\psi_{N}$ to $\Psi_{N}$ on $\R/\Z$:
$$
\Psi_{N}(\gt)
:=
\sum_{m\in\Z}\psi_{N}(\gt+m)
,
$$
and put each such spike at a major arc:
\be
\label{eq:PsiQNis}
\Psi_{\cQ,N}(\gt)
:=
\sum_{q<\cQ}
\sum_{(a,q)=1}
\Psi_{N}\left(\gt-\frac aq\right)
.
\ee
Note that the support of $\Psi_{\cQ,N}$ is $\fM_{\cQ}$.

As in \eqref{eq:RNd}, write
 the representation number
$$
R_{N}(n):=\hat S_{N}(n)=\int_{0}^{1}S_{N}(\gt)e(-n\gt)d\gt,
$$
and decompose it into a (smoothed) major arcs contribution and an error
\be\label{eq:RNdecomp}
R_{N}(n)=\cM_{N}(n)+\cE_{N}(n),
\ee
where
\be
\label{eq:MNdef}
\cM_{N}(n)
:=
\int_{0}^{1}
\Psi_{\cQ,N}(\gt)
S_{N}(\gt)
e(-n\gt)
d\gt
,
\ee
and
\be
\label{eq:ENdef}
\cE_{N}(n)
:=
\int_{0}^{1}
(1-\Psi_{\cQ,N}(\gt))
S_{N}(\gt)
e(-n\gt)
d\gt
.
\ee

The ultimate goal of this section is to prove the following

\begin{thm}
For 
$
\gba
 N\le n< \gbb
 N
,
$
\be\label{eq:cMNbnd}
\cM_{N}(n)
\gg
{1\over \log\log N}
{|\gW_{N}|\over N}
.
\ee
\end{thm}

\pf
Fix 
$
\gba
 N\le n< \gbb
 N
.
$
Starting with \eqref{eq:MNdef}, insert \eqref{eq:PsiQNis} and \eqref{eq:SNsplit} (recall $\supp\Psi_{\cQ,N}\subset\fM_{\cQ}$),
and make the change of variables $\gb=\gt-a/q$.
\bea
\label{eq:MNsplit}
\cM_{N}(n)
&=&
\sum_{q<\cQ}
\sum_{(a,q)=1}
\nu_{q}(a)
e\left(-n\frac aq\right)
\\
\nonumber
&&
\quad
\times
\int_
{0}^{1}
\Psi_{N}
(\gb)\
\vp_{N}
(\gb
)\
e
(-n
\gb
)\
d\gb
\\
\nonumber
&&
\qquad
+
O
\left(
\cQ\cQ
{\cQ\over N}
{|\gW_{N}|}
\cQ^{-4}
\right)
,
\eea
where we used \eqref{eq:vpSize}.

Note that $\cM_{N}$ has already split (up to acceptable error) into the product of the singular series
\be\label{eq:singSer}
\fS_{\cQ}(n)
:=
\sum_{q<\cQ}
\sum_{(a,q)=1}
\nu_{q}(a)
e\left(-n\frac aq\right)
,
\ee
and 
the singular integral
\bea
\nonumber
\Pi_{N}(n)
&:=&
\int_0^{1}
\Psi_{N}(\gb)\
\vp_{N}(\gb)\
e(-n\gb)\
d\gb
=
\sum_{m\in\Z}
\hat
\Psi_{N}
(n-m)
\hat\vp_{N}
(m)
\\
\label{eq:singInt}
&=&
\frac\cQ N
\sum_{m\in\Z}
\hat
\psi
\left(\frac\cQ N(n-m)\right)
\hat\vp_{N}
(m)
.
\eea

First we sketch an analysis of the singular series, which is standard. Insert \eqref{eq:nuQdef} into \eqref{eq:singSer}:
\beann
\fS_{\cQ}(n)
&=&
\sum_{q<\cQ}
\frac1{|\SL_{2}(q)|}
\sum_{\g\in\SL_{2}(q)}
c_{q}\left(\<\g e_{2},e_{2}\>-n\right)
,
\eeann
where $c_{q}$ is the classical Ramanujan sum
$$
c_{q}(m)=\sum_{(a,q)=1}e(am/q)
.
$$
Recall that $c_{q}$ is multiplicative in $q$, and that $c_{q}(m)=\mu(q)$ if $(m,q)=1$ (here $\mu$ is the M\"obius function).
Hence we may extend the range of the sum $q<\cQ$  to $q<\infty$ with a negligible error, obtaining a sum which factors into an Euler product. At each place, the contribution from  prime powers is negligible. We are left to analyze
\bea
\nonumber
\fS_{\cQ}(n)
&\gg&
\fS(n)
\gg
\prod_{p}
\left(
1+
\frac1{|\SL_{2}(p)|}
\sum_{\g\in\SL_{2}(p)}
c_{p}\left(\<\g e_{2},e_{2}\>-n\right)
\right)
\\
\label{eq:singSerBnd}
&=&
\prod_{p\nmid n}
\left(
1+
{1\over p^{2}-1}
\right)
\prod_{p\mid n}
\left(
1-
{1\over p+1}
\right)
\gg
{1\over \log\log n}
.
\eea


Returning to \eqref{eq:singInt}, we now analyze the singular integral. By positivity and using \eqref{eq:psiHatIs} that $\hat\psi(y)>2/5$ for $|y|<1/2$, we have
$$
\Pi_{N}(n)=
\frac\cQ N
\sum_{m\in\Z}
\hat
\psi
\left(\frac\cQ N(n-m)\right)
\hat\vp_{N}
(m)
\ge
\frac25
\frac\cQ N
\sum_{|m-n|<N/(2\cQ)}
\hat\vp_{N}
(m)
.
$$
For $N$ (and hence $\cQ$) sufficiently large, the ranges $n/N\in[\gba,\gbb]$ and $|m-n|<N/(2\cQ)$ force $m/N\in[\gbaa,\gbbb]$, so \eqref{eq:vpBndBelow} applies, giving
\be\label{eq:PiNbnd}
\Pi_{N}(n)
\gg
\frac\cQ N
{N\over2\cQ}
{|\gW_{N}|\over N}
\gg
{|\gW_{N}|\over N}
.
\ee
Inserting \eqref{eq:PiNbnd} and \eqref{eq:singSerBnd} into \eqref{eq:MNsplit} gives \eqref{eq:cMNbnd}, as claimed.
\epf


\section{Minor Arcs Analysis I}\label{sec:Bnd1}

We keep all the notation from the previous section.
Having dealt with the main term \eqref{eq:MNdef}, we are now tasked with estimating the error $\cE_{N}$ in \eqref{eq:ENdef}. As discussed in \eqref{eq:SNN4gd}--\eqref{eq:WQKis}, the key goal is to 
estimate
$$
\sum_{n\in\Z}|\cE_{N}(n)|^{2}=\int_{0}^{1}
|1-\Psi_{\cQ,N}(\gt)|^{2}
|S_{N}(\gt)|^{2}
d\gt
,
$$
where we applied Parseval's formula. 
For $\gt$ outside the major arcs, $\Psi_{\cQ,N}$ vanishes, and we decompose the above integral into regions
$$
W_{Q,K}:=
\left\{
\gt=\frac aq+\gb :
 \foh Q\le 
 q< Q, (a,q)=1,
  \frac K{2N} \le 
  |\gb| < \frac KN
\right\}
.
$$
Here the parameters $Q$ and $K$ range dyadically in 
\be\label{eq:QKrange}
Q<N^{1/2},\qquad 
K<N^{1/2}/Q.
\ee 
If $K=O(1)$, we replace the condition 
$  \foh K/N \le 
  |\gb| <  K/N
$
in $W_{Q,K}$ by just
$
  |\gb| <  K/N,
$
and any appearances of $K$ should be replaced by $1$.

In this section, we give two bounds for $S_{N}(\gt)$, similar to Theorems 5.1 and 6.1 of \cite{BourgainKontorovich2010}.
These will suffice as long as $Q$ or $K$ is large.

\subsection{The bound for $K$ large}\label{sec:Klarge}\

We will first bound $\int_{W_{Q,K}}|S_{N}(\gt)|^{2}d\gt$ by pulling out the largest value of the integrand and multiplying by the measure of the domain, which $\ll Q^{2}\frac KN$. To get the desired bound of $N^{4\gd-1}$ (see \eqref{eq:SNN4gd}), we need to bound the sup norm of $S_{N}$ on $W_{Q,K}$ by a bit less than $N^{2\gd}/(K^{1/2}Q)$. We will win by an extra $K^{1/2}$.

\begin{prop}
Let $N$, $Q$, $K$ be as above, and write $\gt=\frac aq+\gb\in W_{Q,K}$.
Then
\be\label{eq:minBnd1}
|S_N(\gt)|\ll 
N^{2\gd}
\left( {N^{1-\gd}
\over KQ}\right)
,
\ee
as $N\to\infty$.
\end{prop}


\pf
This is a simplified version of Theorem 5.1 in \cite{BourgainKontorovich2010}; we repeat the arguments.
By  \eqref{eq:gWLdef}, we decompose 
\be\label{eq:4.3}
\gW_{N}=
\Xi_{1}
\left(
\Xi_{2}
\Xi_{3}
\cdots
\Xi_{J}
\right)
=
\Xi_{1}\cdot \gW'
.
\ee
Then by \eqref{eq:normToE2}, 
 \eqref{eq:normToTrace}, 
  \eqref{eq:glToTr}, 
Lemma \ref{lem:5.6}, 
and \eqref{eq:NjToN2j},
we have for $\g\in\Xi_{1}$ and $\gw\in\gW'$ that
\be\label{eq:4.4}
|{}^{t}\g e_{2}|
,
|\gw e_{2}|
< 
50 N^{1/2}
.
\ee
Note also from \eqref{eq:XiSize} that
\be\label{eq:Xi1Bnd}
\#\Xi_{1}
,
\# \gW'\ll N^{\gd}
.
 \ee

Then we 
can rewrite $S_{N}(\gt)$ as
\be\label{eq:4.6}
S_{N}(\gt)
=
\sum_{x\in\Z^{2}}
\sum_{y\in\Z^{2}}
\mu(x) \nu (y)
e(\gt \<x,y\>)
,
\ee
where $\mu$ and $\nu$ are image measures in $\Z^{2}$ 
defined by
$$
\mu(x):=\sum_{\g\in\Xi_{1}}\bo_{\{x={}^{t}\g\cdot 
e_{2}\}},
$$
and similarly
$$
\nu(y):=\sum_{\gw\in\gW'}\bo_{\{y=\gw\cdot 
e_{2}\}}.
$$

The projection $\gw\mapsto\gw\cdot e_{2}$ in $\nu$ is 
1-to-1, since it is well-known that the continued fraction of a rational number, if restricted to have even length, is unique.
 The map $\mu$ is also 1-to-1, since $\cG$ is preserved under transposition, ${}^{t}g\in\cG$ for $g\in\cG$ (since its generators \eqref{eq:cGgens} are fixed by transposition).
Hence we have
\be\label{eq:munuInf}
\|\mu\|_{\infty}\le 1,
\qquad
\|\nu\|_{\infty}\le 1
.
\ee

Note that
 for any $y,y'\in\supp\nu$,
 we have  from \eqref{eq:4.4} that $|y-y'|<100N^{1/2}$. Decompose $\nu$ into $100000$ blocks $\nu=\sum_{\ga}\nu^{(\ga)}$ so that for each $\ga$ and any $y,y'\in\supp\nu^{(\ga)}$, 
 \be\label{eq:yypDiff}
 |y-y'|<\frac1{2}N^{1/2}.
 \ee
Write $|S_{N}(\gt)|\le \sum_{\ga}|S_{N}^{(\ga)}(\gt)|$, where
$$
S_{N}^{(\ga)}(\gt):=\sum_{x}\sum_{y}\mu(x)\nu^{(\ga)}(y)e(\gt\<x,y\>).
$$
We will bound each such $S_{N}^{(\ga)}$ independently of $\ga$, so we drop the superscripts $\ga$. 

Let  $\gU:\R^{2}\to \R_{+}$ be a smooth test function which exceeds $1$ on the square $[-1,1]\times[-1,1]$, and has Fourier transform supported in a ball of radius $1$ about the origin.
Apply Cauchy-Schwarz in the $x$ variable, insert $\gU$,
 and open the squares:
$$
|S_N(\gt)|\ll \left(\sum_{x}\mu^{2}(x)\right)^{1/2}
\left(
\sum_{x}\gU\left(\frac x{50N^{1/2}}\right)
\sum_{y}\nu(y)
\sum_{y'}\nu(y')
e(\<x,y-y'\>\gt)
\right)^{1/2}
.
$$
The first parentheses contribute $N^{\gd/2}$ by \eqref{eq:Xi1Bnd} and \eqref{eq:munuInf}.
To the last  sum on $x$ apply Poisson summation, recalling the support of $\hat\gU$:
\be\label{eq:SNPoisson}
|S_N(\gt)|
\ll 
N^{\gd/2}
\left(
\sum_{y}\nu(y)
\sum_{y'}\nu(y')\
N\
\bo_{\{\|(y-y')\gt\|<\frac1{50N^{1/2}}\}}
\right)^{1/2}
.
\ee
Here $\|\cdot\|$ is the distance to the nearest lattice point in $\Z^{2}$. For
such $y,y',\gt$, we have
$$
\|(y-y')\frac aq\|
\le
\|(y-y')\gt\|
 + |y-y'||\gb|
<
 \frac1{50N^{1/2}}
+
\foh N^{1/2}\frac KN
<
\frac1Q
,
$$
where we used \eqref{eq:yypDiff} and \eqref{eq:QKrange}.
Then
$q<Q$
forces
$\|(y-y')\frac aq\|=0$, or
$$
y\equiv y'(q).
$$
This being the case, we now have
$$
\frac1{50N^{1/2}}
>
\|(y-y')\gt\|
=
|(y-y')\gb|
,
$$
that is,
$$
|y-y'|
\ll
\frac{N^{1/2}}K
.
$$
In summary, we have
\beann
|S_N(\gt)|
&\ll&
N^{(\gd+1)/2}
\left(
\sum_{y}\nu(y)
\sum_{y'}
\bo_{\left\{ { y\equiv y'(q) \atop |y-y'|\ll \frac{N^{1/2}}K}\right\}}
\right)^{1/2}
,
\eeann
where we used \eqref{eq:munuInf}. 
Using $Q<\frac{N^{1/2}}K$ and the crudest bound on the $y'$ sum gives
\beann
|S_N(\gt)|
&\ll&
N^{(\gd+1)/2}
\left(
\sum_{y}\nu(y)
\left({N^{1/2}\over QK}\right)^{2}
\right)^{1/2}
\ll
{N^{\gd+1}
\over QK}
,
\eeann
as claimed.
\epf

The bound \eqref{eq:minBnd1} is already conclusive if  
$
K
$
is a bit larger than $N^{2(1-\gd)}
.
$
\begin{thm}\label{thm:minBnd1}
Assume 
$Q<N^{1/2}$ and  $K<N^{1/2}/Q$. 
Then
\be\label{eq:minBnd1'}
\int_{W_{Q,K}}|S_N(\gt)|^{2}d\gt \ll 
{
(\#\gW_{N})^{2}
\over N}
\left[{N^{2(1-\gd)+4\Bpow} 
\over K}\right]
.
\ee
\end{thm}
\pf
We bound trivially using \eqref{eq:minBnd1}:
$$
\int_{W_{Q,K}}|S_N(\gt)|^{2}d\gt
\ll
\frac KN Q^{2}
\left(
{N^{\gd+1}
\over QK}
\right)^{2}
\ll
N^{4\gd-1}\left[{N^{2(1-\gd)}\over K}\right]
,
$$
and the claim follows from \eqref{eq:gWLsize}, on crudely using $\gd<1$.
\epf

\subsection{Another Bilinear Forms Estimate}\label{sec:Bilin2}\

Next we introduce the cross-section of $W_{Q,K}$ for fixed $\gb$:
$$
P_{Q,\gb}:=\left\{\gt=\frac aq+\gb: \foh Q\le q< Q, (a,q)=1\right\}.
$$
We will bound using \eqref{eq:minBnd1}, giving essentially
$$
\int_{W_{Q,K}}|S_{N}|^{2}
\ll
\sup|S_{N}|\frac KN\sup_{\gb}\sum_{P_{Q,\gb}}|S_{N}|
\ll
{N^{2\gd+}\over KQ}\frac KN\sup_{\gb}\sum_{P_{Q,\gb}}|S_{N}|
.
$$
The trivial bound on $\sum_{P_{Q,\gb}}|S_{N}|$ is of course $N^{2\gd}Q^{2}$, so we need to save
a little more than a power of $Q$ to get our target bound of less than $N^{4\gd-1}$. 
This is
achieved by exploiting the extra structure in the $a$ and $q$ sums, as follows.

\begin{prop}\label{prop:PQgb}
Let the notation be as above.
 Then
 for all $\vep>0$,
\be\label{eq:minBnd2}
\sum_{\gt\in P_{Q,\gb}}|S_N(\gt)| \ll_{\vep}
N^{2\gd} Q^{2}
N^{1-\gd+\vep}
\left[
{1
\over   Q^{3/2}}
+
{
1
\over
Q
N^{1/8}
}
\right]
.
\ee
\end{prop}


\pf
The proof is 
nearly
 identical to that of Theorem 6.1 in \cite{BourgainKontorovich2010}, but we reproduce it for the reader's convenience. 
 We again use \eqref{eq:gWLdef} to decompose $\gW_{N}$ into pieces, now grouping by
 $$
 \gW_{N}=
 \left(\Xi_{1}\Xi_{2}\right)
 \left(\Xi_{3}\cdots\Xi_{J}\right)
 =
 \gW'\cdot\gW''
 .
 $$
 As before, we have for $\g\in\gW'$ and $\gw\in\gW''$ that
 \be\label{eq:4.4p}
 |{}^{t}\g e_{2}|<300N^{3/4},
\qquad
 \text{ and }
\qquad
 |\gw e_{2}|<2000 N^{1/4}
 .
 \ee
 Also from \eqref{eq:XiSize}, we have
 \be\label{eq:gWpBnd}
\# \gW'\ll N^{3\gd/2}
\qquad
 \text{ and }
\qquad
\# \gW''\ll N^{\gd/2}
.
 \ee
 Again we define the measures $\mu$ and $\nu$ on $\Z^{2}$ by
 $$
 \mu(x):=\sum_{\g\in\gW'}\bo_{\{x={}^{t}\g e_{2}\}}
,
 $$
 $$
 \nu(y):=\sum_{\gw\in\gW''}\bo_{\{y=\gw e_{2}\}}
,
 $$
 with $\mu,\nu\le1$. For any two elements $y,y'$ in the support of $\nu$, we have $|y-y'|<4000N^{1/4}$. Hence we again decompose $\nu$ into $O(1)$ pieces, $\nu=\sum_{\ga}\nu^{(\ga)}$, so as to make the difference 
 \be\label{eq:yypDiffNew}
|y-y'|<\frac1{10000} N^{1/4},
 \ee
for $y,y'$ in the support of $\nu^{(\ga)}$. Writing
$$
S_{N}^{(\ga)}(\gt)=\sum_{x}\sum_{y}\mu(x)\nu^{(\ga)}(y)e(\gt\<x,y\>),
$$
and dropping the superscripts $\ga$,
we proceed to bound
\beann
\sum_{\gt\in P_{Q,\gb}}|S_N(\gt)| 
&=&
\sum_{q\asymp Q}
\sum_{(a,q)=1}
\gz(\gt)
S_N(\gt)
\\
&=&
\sum_{q\asymp Q}
\sum_{(a,q)=1}
\gz(\gt)
\sum_{x}
\sum_{y}
\mu(x)
\nu(y)
e(\gt\<x,y\>)
,
\eeann
where $\gz$ has modulus $1$. 
Recall the bump function $\gU$ which is at least one on $[-1,1]^{2}$; assume now that its Fourier transform is supported in a ball of radius $1/40$ about the origin.
Apply Cauchy-Schwarz in the $x$ sum and \eqref{eq:gWpBnd}, insert the  function $\gU$, reverse orders, and apply Poisson summation:
\bea
\nonumber
\sum_{\gt\in P_{Q,\gb}}|S_N(\gt)| 
&\ll&
N^{3\gd/4}
\left(
\sum_{x}
\gU\left({x\over 300N^{3/4}}\right)
\left|
\sum_{q\asymp Q}
\sum_{(a,q)=1}
\gz(\gt)
\sum_{y}
\nu(y)
e(\gt\<x,y\>)
\right|^{2}
\right)^{1/2}
\\
\label{eq:PbetaPoisson}
&\ll&
N^{3(\gd+1)/4}\
\cX
^{1/2}
,
\eea
where
\be\label{eq:cXis}
\cX=\cX_{
Q,\gb}:=
\sum_{q}
\sum_{q'}
\sum_{a}
\sum_{a'}
\sum_{y}
\sum_{y'}
\nu(y)
\nu(y')
\bo_{\left\{\|y\gt-y'\gt'\|<\frac1{12000N^{3/4}}\right\}}
.
\ee
Here $\gt'=\frac{a'}{q'}+\gb$; note that $\gb$ is the same for $\gt$ and $\gt'$.

Write $y=(y_{1},y_{2})$ and the same with $y'$. Consider the innermost condition in \eqref{eq:cXis}: 
\be\label{eq:y1y2gt}
\|y_{1}\gt-y_{1}'\gt'\|
,
\|y_{2}\gt-y_{2}'\gt'\|<\frac1{12000N^{3/4}}
.
\ee
Recall that $y=\g e_{2}$ for some 
(non-identity) $\g\in\G$, and the same for $y'$; hence we have
$$
y_{1}y_{2}
y'_{1}y'_{2}
\neq0
.
$$
Also note using \eqref{eq:y1y2gt}, \eqref{eq:yypDiffNew} and $|\gb|<K/N<1/(N^{1/2}Q)$ that 
\be\label{eq:ygtmypgtp}
\left\|y_{1}\frac aq-y_{1}'\frac{a'}{q'}\right\|
\le 
\|y_{1}\gt-y_{1}'\gt'\|
+
|(y_{1}-y_{1}')\gb|
<
\frac1{12000N^{3/4}}
+
\frac{N^{1/4}}{10000N^{1/2}Q}
,
\ee
and similarly with $y_{2},y_{2}'$.




Let $Y:=\mattwo{y_{1}}{y_{1}'}{y_{2}}{y_{2}'}$, so that 
\be\label{eq:cYdef}
\cY:=\det(Y)=y_{1}y_{2}'-y_{1}'y_{2}.
\ee
Observe then by \eqref{eq:ygtmypgtp}, \eqref{eq:4.4p}, and $Q<N^{1/2}$ that
\beann
\left\|
\cY
\frac aq
\right\|
&\le&
\left\|
y_{2}'
\left(y_{1}
\frac aq
-
y_{1}'
\frac {a'}{q'}
\right)
\right\|
+
\left\|
y_{1}'
\left(
y_{2}'
\frac {a'}{q'}
-
y_{2}\frac aq
\right)
\right\|
\\
&<&
2000N^{1/4}
\left(
\frac1{12000N^{3/4}}
+
\frac{N^{1/4}}{10000N^{1/2}Q}
\right)
\times2
\\
&<&
\frac1Q
.
\eeann
Of course this forces
$
\cY
\equiv 0 (\mod q)$. 
The same argument gives
$
\cY
\equiv 0 (\mod q')$, and hence we have
\be\label{eq:cYmodq}
\cY
\equiv 0 (\mod \fq)
,
\ee
where $\foh Q\le \fq<Q^{2}$ is the least common multiple of $q$ and $q'$.

Decompose $\cX$ in \eqref{eq:cXis} as $\cX=\cX_{1}+\cX_{2}$ according to
whether $
\cY
=0$ or not; we handle the two contributions separately.
We will prove
\begin{lem}\label{lem:X1bnd}
For any $\vep>0$,
$$
\cX_{1}
\ll_{\vep}
N^{\gd/2+\vep}
Q^{4}
\left[
\frac1{N^{3/4}}
+
Q^{-2}
\right]
.
$$
\end{lem}
and
\begin{lem}\label{lem:X2bnd}
For any $\vep>0$,
$$
\cX_{2}
\ll_{\vep}
N^{\gd+\vep}
Q
.
$$
\end{lem}

We momentarily postpone the proofs of
these two Lemmata,
first using them to  finish the proof of Proposition \ref{prop:PQgb}.
Returning to \eqref{eq:PbetaPoisson}, we have
\beann
\sum_{\gt\in P_{Q,\gb}}
|S_{N}(\gt)|
&\ll_{\vep}&
N^{3(\gd+1)/4+\vep}
\left[
N^{\gd/2}
Q^{4}
\left(
\frac1{N^{3/4}}
+
Q^{-2}
\right)
+
N^{\gd}
Q
\right]^{1/2}
,
\eeann
from which the claim follows using $Q<N^{1/2}$.
\epf

Now we establish the Lemmata separately.

\subsubsection{Bounding $\cX_{2}$: the case $
\cY
\neq0$}\

\pf[Proof of Lemma \ref{lem:X2bnd}]
Note from \eqref{eq:4.4p}, \eqref{eq:yypDiffNew}, and \eqref{eq:cYdef} that 
$$
|\cY|
\le
|y_{1}(y_{2}'-y_{2})|
+
|(y_{1}-y_{1}')y_{2}|
<
2000N^{1/4}
{1\over 10000}N^{1/4}
\times2
<N^{1/2}
.
$$
Since $\fq\mid 
\cY
$ and $\cY\neq0$, we have 
$$
\fq\le \min(Q^{2}, 
N^{1/2})\le QN^{1/4}
.
$$
Then \eqref{eq:ygtmypgtp} and $Q<N^{1/2}$ forces 
\be\label{eq:yypmod1}
y_{1}\frac aq-y_{1}'\frac{a'}{q'}
\equiv 
0(\mod 1),
\ee
and the same holds for $y_{2},y_{2}'$. Let $\tilde q:=(q,q')$ and $q=q_{1}\tilde q$, $q'=q_{1}'\tilde q$ so that $\fq=q_{1}q_{1}'\tilde q$.
Then \eqref{eq:yypmod1} becomes
$$
y_{1}
a
q_{1}'
\equiv 
y_{1}'
{a'}
q_{1}
(\mod
\fq
)
,
$$
and the same for $y_{2},y_{2}'$.
Recall $a$ and $q$ are coprime, as are $a'$ and $q'$.
%
It then follows that $q_{1}\mid y_{1}$, and similarly, $q_{1}\mid y_{2}$. But since $y$ is a visual vector, $(y_{1},y_{2})=1$, forcing $q_{1}=1$. 
The same argument applies to $q_{1}'$, so we have $q=q'=\fq$.
Then \eqref{eq:yypmod1} now reads
\be\label{eq:y1apdet}
y_{1}a\equiv y_{1}'a'
(\mod q)
,
\ee
and similarly for $y_{2},y_{2}'$.

Hence, once we fix $y,y'\in\gW'' e_{2}$, the value of $\cY$ is determined, and $q\mid \cY$ leaves at most $N^{\vep}$ choices for $q$.
Then there are at most $Q$ choices for $a$, from which $a'$ is determined by \eqref{eq:y1apdet} (again using that $y$ and $y'$ are visual vectors).

Then using \eqref{eq:gWpBnd}, $\cX_{2}$ is bounded by
\bea
\label{eq:cX2bnd}
\cX_{2}
&\ll&
\sum_{y}
\nu(y)
\sum_{y'}
\nu(y')
\sum_{q\mid\cY\atop \foh Q\le q<Q}
\sum_{a(\mod q)}
1
\\
\nonumber
&\ll_{\vep}&
\left(N^{\gd/2}\right)^{2}
N^{\vep}
Q
,
\eea
as claimed.
\epf

\subsubsection{Bounding $\cX_{1}$: the case $
\cY
=0$}\

\pf[Proof of Lemma \ref{lem:X1bnd}]
The condition $\cY=0$ implies $y_{1}/y_{2}=y_{1}'/y_{2}'$. Recall that rationals have unique continued fraction expansions (of even length), and thus $y=y'$.
The bottom line savings from this fact is at most $N^{1/4}$, whereas we need to save a bit more than $Q$, which can be as large as $N^{1/2}$.

Let $N':=\frac1{12000}N^{3/4}$.
The condition \eqref{eq:y1y2gt} then becomes
\be\label{eq:y1gtnew}
\left\|y_{1}\left(\frac aq-\frac{a'}{q'}\right)
\right\|
<
\frac 1{N'}
.
\ee
Let $(y_{1},q)=v$ and $(y_{1},q')=v'$ with $q=vr$. Assume without loss of generality that $v\le v'$.
Fix $y$ (for which there are $N^{\gd/2}$ choices) and $v,v'\mid y_{1}$ (at most $N^{\vep}$ choices). There are $\ll Q/v'$ choices for $q'\equiv0(\mod v')$, and then $\ll Q$ choices for $(a',q')=1$. Write $\psi$ for $y_{1}a'/q'(\mod 1)$, which is now fixed, and write $y_{1}=vz$ with $(z,r)=1$.
Then \eqref{eq:y1gtnew} becomes
$$
\left\|
z\frac ar-\psi
\right\|
<
\frac1{N'}
.
$$
Let $\cU_{z}$ be the set of possible fractions ${za\over r} (\mod 1)$ as $r$ varies in 
$
{Q/ (2v)}\le r<{Q/ v}
,
$ 
and $a$ ranges up to $Q$ subject to $(a,vr)=1$.
Note that distinct points $u\in\cU_{z}$ are separated by a distance of at least $
v^{2}/Q^{2}$. Hence the size of the intersection of $\cU_{z}$ with the interval 
$$
\left[\psi-\frac1{N'},\psi+\frac1{N'}\right]
$$ 
contains at most $\frac{Q^{2}}{v^{2}N'}+1$ points. Once $u=f/r\in\cU_{z}$ is determined, so is its denominator, that is, $r$ is determined. Also $a(\mod r)$ is determined (to be $f$), hence  $a(\mod q)$ has $v$ possible values (recall $q=rv$).

In summary, we use \eqref{eq:gWpBnd} again to bound $\cX_{1}$ by:
\beann
\cX_{1}
&\ll&
\sum_{y}
\nu(y)
\sum_{{v, 
v'\mid y_{1}}\atop v\le v'}
\sum_{q'\equiv 0(\mod v')}
\sum_{(a',q')=1}
\sum_{f/r\in\cU_{z}\cap[\psi-\frac1{N'},\psi+\frac1{N'}]\atop q=rv}
\sum_{a<q\atop a\equiv f(\mod r)}
1
\\
&\ll&
\sum_{y}
\nu(y)
\sum_{{v, 
v'\mid y_{1}}\atop v\le v'}
\frac Q{v'}
Q
\left(
\frac{Q^{2}}{v^{2}N'}+1
\right)
v
\\
&\ll_{\vep}&
N^{\gd/2}
N^{\vep}
Q^{2}
\left(
\frac{Q^{2}}{N^{3/4}}+1
\right)
,
\eeann
as claimed.
\epf

With the Lemmata established, we have completed the proof of Propostion \ref{prop:PQgb}.

\subsection{The bound for $Q$ large}\

Lastly, we input this bound to get another bound on the main integral, one which is favorable as long as $Q$ is a bit bigger than $N^{4(1-\gd)}$. 

\begin{thm}\label{thm:minBnd2}
Assume that 
$Q<N^{1/2}$ and $KQ<N^{1/2}$.
 Then 
\be\label{eq:minBnd2'}
\int_{W_{Q,K}}|S_N(\gt)|^{2}d\gt 
\ll
{(\#\gW_{N})^{2}\over N}
N^{2(1-\gd)
} 
N^{4
\Bpow}
 \left(
{1 \over  Q^{1/2}}
+
{1 \over  N^{1/8}}
\right)
.
\ee
\end{thm}

\pf
Write
\beann
\int_{W_{Q,K}}|S_N(\gt)|^{2}d\gt
&\ll&
\sup_{\gt\in W_{Q,K}}|S_N(\gt)|
\cdot
\frac KN
\sup_{|\gb|\asymp \frac KN}
\sum_{\gt\in P_{Q,\gb}}|S_N(\gt)| 
\\
&\ll_{\vep}&
N^{2\gd}\left({N^{1-\gd}\over KQ}\right)
\cdot
\frac KN
\left(
N^{2\gd}
Q^{2}
N^{1-\gd+\vep}
\left[
\frac1{Q^{3/2}}
+
\frac1{QN^{1/8}}
\right]
\right)
\\
&\ll&
N^{4\gd-1}
N^{2(1-\gd)+\vep}
\left({
1
\over Q^{1/2}
}
+
{
1\over
N^{1/8}
}
\right)
,
\eeann
where we used \eqref{eq:minBnd1} and  \eqref{eq:minBnd2}. The claim follows from \eqref{eq:gWLsize}, again crudely using $\gd<1$.
\epf

It remains to handle the regions when both $K$ and $Q$ are very small, less than $N^{\vep}$ for $\vep$ near zero. 



\section{Minor Arcs Analysis II}\label{sec:Bnd2}

We now push the methods of the previous section down to the level of $Q$ and $K$ being of constant size. We again do this in two stages. But first we record the following counting bound. 

\begin{lem}\label{lem:BKS}
For $(qK)^{\gdquart}<Y<X$, and visual vectors $\eta,\eta'\in\Z^{2}$ (meaning their coordinates are coprime) with  $|\eta|\asymp X/Y$ and $|\eta'|\asymp Y$, 
$$
\#
\left\{\g\in\SL_{2}(\Z): \|\g\|\asymp Y, |\g\eta-\eta'|< \frac X{YK}, \text{ and } \g\eta\equiv\eta'(\mod q)\right\}
\ll
\left({Y\over qK}\right)^{2}
.
$$
The implied constant is absolute, depending on the implied constants above.
\end{lem}
\pf[Sketch of proof]
Write $G(\Z)=\SL_{2}(\Z)$ and let $G_{\eta}(q)$ be the stabilizer of $\eta$ mod $q$:
$$
G_{\eta}(q):=\{\g\in G(\Z):\g\eta\equiv\eta(q)\}.
$$
Then $G(\Z)\cong(G(\Z)/G_{\eta}(q))\times G_{\eta}(q)$.
Let $R=R_{Y,K}$ denote the region 
$$
R:=
\{g\in\SL_{2}(\R):
\|g\|\asymp Y, |g\eta-\eta'|<  X/(YK)
\}
.
$$
The methods in \cite{GoodBook} (see also \cite{BourgainKontorovichSarnak2010})
 give an estimate of the form
\beann
\sum_{\g\in G(\Z)}\bo_{\{\g\in R\}}\bo_{\g\eta\equiv\eta'(\mod q)}
&=&
\sum_{\gw\in G(\Z)/G_{\eta}(q)}
\bo_{\gw\eta\equiv\eta'(\mod q)}
\sum_{\g'\in G_{\eta}(q)}\bo_{\{\gw\g'\in R\}}
\\
&\ll_{\vep}&
\sum_{\gw\in G(\Z)/G_{\eta}(q)}
\bo_{\gw\eta\equiv\eta'(\mod q)}
\left(
\left({Y\over qK}\right)^{2}
+
Y^{2\gT+\vep}
\right)
\\
&\ll&
\left({Y\over qK}\right)^{2}
+
Y^{2\gT+\vep}
,
\eeann
where $\gT=1/2+7/64$ is the best known bound towards the Ramanujan conjectures \cite{KimSarnak2003}. (We apply the argument to a  smoothed sum.)
The first term
dominates as long as $(qK)^{2}<Y^{25/32-\vep}$, and the claim follows using $64/25+\vep<65/25=\gdquart$. 
\epf
\begin{rmk}\label{rmk:bestDel}
Recall that we
are not interested here in
 optimizing the final value of
  $\gd_{0}$
  in Theorem \ref{thm:mainN},
 so allow ourselves to be a bit crude in the above for the sake of exposition. 
\end{rmk}

\subsection{The bound for $K$ at least a small power of $Q$}\

We return to the approach of \S\ref{sec:Klarge}, that is just bounding the sup norm and needing to win more than $K^{1/2}Q$ off the trivial bound. Now we use the fact that $KQ$ is quite small to beat the trivial bound by $(KQ)^{1-\vep}$ with $\vep$ small (depending on the distance from $\gd$ to $1$).
Then we will have, roughly 
$$
\int_{W_{Q,K}}|S_{N}|^{2}
\ll 
Q^{2}\frac KN\left({\#\gW_{N}\over (KQ)^{1-\vep}}\right)^{2}
=
{(\#\gW_{N})^{2}\over N}
\left({Q^{\vep}\over K^{1-\vep}}\right)
,
$$
which is a savings as long as $K$ is at least a small power of $Q$. Note now for $K$ and $Q$ small that we must be careful with the loss in the size of $\gW_{N}$ in the lower bound \eqref{eq:gWLsize}.
We make all of this precise below.

\begin{prop}
Assume $\gt\in W_{Q,K}$ with
\be\label{eq:KQsm}
1\ll KQ<N^{\gdover}
.
\ee 
 Then
\be\label{eq:Sbnd}
|S_N(\gt)|\ll
\#\gW_{N}
\left(
{
e^{c(\log\log (KQ))^{2}}
\over (KQ)^{1-(1-\gd)\gdone}}
\right)
.
\ee
\end{prop}

\pf
Recalling \eqref{eq:NjToN2j} that $N_{j}\asymp N^{1/2^{j}}$, we can find a $1\le j\le J$ so that
\be\label{eq:KQM}
\frac1{100}(QK)^{\gdquart}<N_{j}<(QK)^{\gdhalf},
\ee
say. Here we used \eqref{eq:KQsm} that $(QK)^{\gdhalf}<N^{1/2}$. 

 Define the sets
\bea
\label{eq:123decomp}
\gW^{(1)}
&:=&
\Xi_{1}
\Xi_{2}
\cdots
\Xi_{j-1}
,
\\
\nonumber
\gW^{(2)}
&:=&
\Xi_{j}
,
\\
\nonumber
\gW^{(3)}
&:=&
\Xi_{j+1}
\Xi_{j+2}
\cdots
\Xi_{J}
.
\eea
Hence for $g_{i}\in\gW^{(i)}$,
\bea
\label{eq:5.3}
\gl(g_{3}) 
&\sim&
N_{j+1}
N_{j+2}
\cdots
N_{J}
=:M
,
\\
\label{eq:5.4}
\gl(g_{2})
&\sim&
N_{j}
,
\\
\label{eq:5.5}
\gl(g_{1})
&\asymp&
{N\over M\ N_{j}}
.
\eea
Note that
\be\label{eq:MtoNj}
{ N_{j}\over \log N_{j}}\ll 
M\ll N_{j}\log N_{j},
\ee 
and that
from \eqref{eq:GWLsize2} and \eqref{eq:XiSize}
we have
\be\label{eq:5.6}
|\gW^{(3)}|
\gg {M^{2\gd}\over e^{c(\log \log M)^{2}}}
,
\qquad
|\gW^{(2)}|
\gg{ (N_{j})^{2\gd}\over (\log N_{j})^{3}}
.
\ee
In the above, we used that $J-j\asymp \log\log M$.

Estimate
\be\label{eq:5.7}
|S_N(\gt)|
\ll
\sum_{g_{1}\in\gW^{(1)}}
\sum_{g_{3}\in\gW^{(3)}}
\left|
\sum_{g_{2}\in\gW^{(2)}}
e(\<g_{3}e_{2},{}^{t}g_{2}{}^{t}g_{1}e_{2}\>\gt)
\right|
.
\ee
Fix $g_{1}$ and set $\eta={}^{t}g_{1}e_{2}$. Note that 
\be\label{eq:etaSize}
|\eta|\asymp \frac N{MN_{j}}.
\ee 
Estimate
as in \eqref{eq:SNPoisson}:
\bea
\label{eq:5.8}
&&
\sum_{g_{3}\in\gW^{(3)}}
\left|
\sum_{g_{2}\in\gW^{(2)}}
e(\<g_{3}e_{2},{}^{t}g_{2}\eta\>\gt)
\right|
\\
\nonumber
&&
\qquad
\ll
(\#\gW^{(3)})^{1/2}\ 
M
\left[
\#
\left\{
(g,g')\in{}^{t}\gW^{(2)}\times{}^{t}\gW^{(2)}
:
{
\|\<(g-g')\eta,e_{1}\>\| \ll \frac1M
\atop
\|\<(g-g')\eta,e_{2}\>\| \ll \frac1M
}
\right\}
\right]^{1/2}
,
\eea
where we extended the sum over $g_{3}$  to  $g_{3}e_{2}\in \{z\in\Z^{2}:|z|\ll M\}$.
Write
\be\label{eq:5.9}
\left\|\<(g-g')\eta,e_{i}\>\frac aq\right\|
=
\left\|
\<(g-g')\eta,e_{i}\>\gt
\right\|
+
\left|
\<(g-g')\eta,e_{i}\>\gb
\right|
,
\ee
where
$$
|
\<(g-g')\eta,e_{i}\>\gb
|
\ll
N_{j}
{N\over M N_{j}}
{K\over N}
=
\frac KM
.
$$
From \eqref{eq:KQM} and \eqref{eq:MtoNj} we clearly have $\frac KM<\frac1Q$, so  \eqref{eq:5.9} forces
\be\label{eq:5.10}
(g-g')\eta\equiv0(q)
\ee
and
\be\label{eq:5.11}
|(g-g')\eta|\ll \frac1{M|\gb|}\ll {N\over KM}
.
\ee
Fix $g'$ and enlarge $g\in{}^{t}\gW^{(2)}$  to $\{g\in\SL_{2}(\Z):\|g\|\ll N_{j}\}$. 
Applying
Lemma \ref{lem:BKS} with $\eta'=g'\eta$, $X=N/M$, and $Y=N_{j}$, 
the $g$ cardinality contributes 
$$
\ll
\left({N_{j}\over KQ}\right)^{2}
.
$$
Thus
we have by \eqref{eq:5.6} and \eqref{eq:MtoNj} that
$$
\eqref{eq:5.8}
\ll
\left(\#\gW^{(3)}\cdot\#\gW^{(2)}\right)^{1/2}
{N_{j}^{2}\over KQ}
\ll
\#\gW^{(3)}\cdot\#\gW^{(2)}
{
(MN_{j})^{1-\gd}
e^{c(\log\log M)^{2}}
(\log N_{j})^{3}
\over KQ}
.
$$
Hence by \eqref{eq:KQM} and \eqref{eq:MtoNj},
$$
\eqref{eq:5.7}
\ll
\#\gW_{N}
{
(KQ)^{(1-\gd)\gdone}
e^{c(\log\log (KQ))^{2}}
\over KQ}
,
$$
as claimed.
\epf

Inserting this bound into the main integral and estimating trivially gives
\begin{thm}\label{thm:minBnd3}
Assuming \eqref{eq:KQsm},
$$
\int_{W_{Q,K}}
|S_N(\gt)|^{2}d\gt
\ll
{(\#\gW_{N})^{2}\over N}
{
Q^{(1-\gd)\gdtwo}
e^{c(\log\log(KQ))^{2}}
\over K^{1-(1-\gd)\gdtwo}}
.
$$
\end{thm}

This bound is conclusive unless 
$
K$
is much less than
\be\label{eq:KllQe}
Q^{\gdtwo(1-\gd)\over 1-\gdtwo(1-\gd)}\approx Q^{\vep}
.
\ee

\subsection{The bound for $K$ even smaller}\

In this last section, we give the final bound for minor arcs, which we apply to the remaining range of $K$ much less than \eqref{eq:KllQe}.
Recall the approach of \S\ref{sec:Bilin2}: we bound $\int_{W_{Q,K}}|S_{N}|^{2}$ by the sup norm times $K/N$ times $\sum_{P_{\gb,Q}}|S_{N}|$. The sup norm has already won almost $KQ$, so we need to win more than a power of $Q$ off of the last summation. We proceed as follows.

\begin{prop}
\label{prop:bnd4}
Recall the cross section $P_{Q,\gb}$ for a fixed $|\gb|\asymp\frac KN$:
$$
P_{Q,\gb}:=\left\{\gt=\frac aq+\gb: q\asymp Q, (a,q)=1\right\}.
$$
Then assuming \eqref{eq:KQsm}, 
we have
\be\label{eq:6.2}
\sum_{\gt\in P_{Q,\gb}}
|S_N(\gt)|
\ll
\#\gW_{N}\
Q^{2}
\left(
{
(KQ)^{(1-\gd)\gdone}
e^{c(\log\log(KQ))^{2}}
\over
Q^{3/2}
}
\right)
.
\ee
\end{prop}

\pf 
This argument is similar to Proposition \ref{prop:PQgb} and we sketch the proof. 
Using the same decomposition \eqref{eq:123decomp}, we follow
 \eqref{eq:PbetaPoisson} and bound the left hand side of \eqref{eq:6.2} by
\bea
\label{eq:6.4}
&\ll& 
\sum_{g_{1}\in\gW^{(1)}}
\left(\#\gW^{(3)}\right)^{1/2}
M
\\
\nonumber
&&
\times
\left[
\#
\left\{
(\gt,\gt',g,g')\in P_{\gb}\times P_{\gb}\times {}^{t}\gW^{(2)}\times {}^{t}\gW^{(2)}
:
\|
(g\gt-g'\gt')\eta
\|
\ll \frac 1M
\right\}
\right]^{1/2}
,
\eea
where $\eta={}^{t}g_{1}e_{2}$.
The innermost condition guarantees $q=q'$ and
$$
a(g\eta)\equiv a'(g'\eta) (\mod q)
,
$$
The number of choices for $g'$ given $g$, $q$, $a$, and $a'$ is
$$
\ll 
\left(
{N_{j}\over Q
}
\right)^{2}
,
$$
hence
\beann
\eqref{eq:6.4}
&\ll&
\#\gW^{(1)}
\left(\#\gW^{(3)}\right)^{1/2}
M
\left[
Q
Q^{2}
\#\gW^{(2)}
\left(
{N_{j}\over Q
}
\right)^{2}
\right]^{1/2}
\\
&\ll&
\#\gW_{N}\
Q^{2}
\left(
{
(KQ)^{(1-\gd)\gdone}
e^{c(\log\log(KQ))^{2}}
\over
Q^{3/2}
}
\right)
,
\eeann
as claimed.
\epf

Using \eqref{eq:6.2} and \eqref{eq:Sbnd}, we now have the bound:
\beann
\int_{W_{Q,K}}|S_N(\gt)|^{2}d\gt
&\ll&
\frac KN
\#\gW_{N}
\left(
{
e^{c(\log\log (KQ))^{2}}
\over (KQ)^{1-(1-\gd)\gdone}}
\right)
\\
&&
\qquad\times
\#\gW_{N}\
Q^{2}
\left(
{
(KQ)^{(1-\gd)\gdone}
e^{c(\log\log(KQ))^{2}}
\over
Q^{3/2}
}
\right)
,
\eeann
from which we immediately have:

\begin{thm}\label{thm:minBnd4}
Assuming \eqref{eq:KQsm},
\beann
\int_{W_{Q,K}}|S_N(\gt)|^{2}d\gt
&\ll&
{(\#\gW_{N})^{2}\over N}
\left(
{(KQ)^{(1-\gd)\gdtwo}
e^{c(\log\log(KQ))^{2}}
\over
Q^{1/2}
}
\right)
.
\eeann
\end{thm}

\section{Proofs of Theorems \ref{thm:mainN} and \ref{thm:primes}}\label{sec:posProp}  

Keeping all the notation of previous sections, we now prove the main minor arcs estimate, analogous to \eqref{eq:SNN4gd}, 
before completing a proof of Theorem \ref{thm:mainN}.

\begin{thm}
Assume
\be\label{eq:gdBnd}
\gd>\gd_{0},
\ee 
with $\gd_{0}$ given by \eqref{eq:gdAllow}.
Then for some $c>0$,
\be\label{eq:ErrBnd}
\sum_{n\in\Z}|\cE_{N}(n)|^{2}
\ll
{|\gW_{N}|^{2}\over N}
\cQ^{-c}.
\ee
\end{thm}
\pf
By Parseval, we have
$$
\sum_{n\in\Z}|\cE_{N}(n)|^{2}=
\int_{0}^{1}
|1-\Psi_{\cQ,N}(\gt)|^{2}|S_{N}(\gt)|^{2}d\gt
=
\int_{\fM_{\cQ}}+\int_{\fm}
,
$$
where we broke the integral into the major arcs $\fM_{\cQ}$ and the complementary minor arcs $\fm=[0,1]\setminus\fM_{\cQ}$.

On the major arcs, note 
from \eqref{eq:psiIs}
that $1-\psi(x)=|x|$ on $[-1,1]$. Then using \eqref{eq:Sbnd} 
with $K\asymp N|\gb|$ gives 
\beann
\int_{\fM_{\cQ}}
&\ll&
\sum_{q<\cQ}
\sum_{(a,q)=1}
\int_{|\gb|<\cQ/N}
\left|
{N\over \cQ}
\gb
\right|^{2}
\left(
|\gW_{N}|
\left(
{
1
\over (N|\gb|\cQ)^{1-c}}
\right)
\right)^{2}
d\gb
\\
&\ll&
{|\gW_{N}|^{2}\over N}
{1\over\cQ^{1-4c}}
.
\eeann
Here $0<c<(1-\gd)52/5<1/4$ by \eqref{eq:gdBnd}, so
 renaming the constant $c>0$, we are done with the major arcs.

Decompose the minor arcs $\fm$ into dyadic regions
$$
\int_{\fm}
|S_{N}(\gt)|^{2}d\gt
\ll
\sum_{ Q <N^{1/2}\atop \text{dyadic}}
\sum_{ K <{N^{1/2}\over Q}\atop \text{dyadic}}
\cI_{Q,K}
,
$$
where
at least one of $Q$ or $K$ exceeds $\cQ$, and
$$
\cI(Q,K):=
\int_{W_{Q,K}}
|S_{N}(\gt)|^{2}d\gt
.
$$

Write 
$Q=N^{\ga}$, $K=N^{\gk}$, with the parameters $(\ga,\gk)$ ranging in 
\be\label{eq:gagkRange}
0\le\ga<1/2\text{ and }0\le\gk<1/2-\ga.
\ee 
It will be convenient to define 
\be\label{eq:etaIs}
\eta:=(1-\gd)\gdtwo
.
\ee
Assume that $1-\gd<\gdofour$, so that $0<\eta<1/2$.
We break the summation into the following four ranges:
\beann
\cR_{1}
&:=&
\{(\ga,\gk):\gk>2(1-\gd)+4\Bpow\}
,
\\
\cR_{2}
&:=&
\{(\ga,\gk):\ga>4(1-\gd)+8\Bpow\}
,
\\
\nonumber
\cR_{3}
&:=&
\{(\ga,\gk): \eta(\ga+\gk)<\gk \text { and } \ga+\gk<\gdover\}
,
\\
\nonumber
\cR_{4}
&:=&
\{(\ga,\gk): \eta(\ga+\gk)<\foh\ga \text { and } \ga+\gk<\gdover\}
.
\eeann

We need to show  
that these four regions cover the entire range \eqref{eq:gagkRange}. 
Using  \eqref{eq:gbIs} and \eqref{eq:gdBnd} with \eqref{eq:gdAllow} guarantees that
the regions $\cR_{1}$ and $\cR_{2}$ certainly cover the range $\ga+\gk\ge \gdover$. 
In the complimentary range, $\cR_{3}$ and $\cR_{4}$ give two regions: the region below the line through the origin with slope $\eta/(1-\eta)$, and the region above the line through the origin with slope $(1/2-\eta)/\eta.$ These two regions overlap when the slopes overlap, that is, when $\eta<1/3$. 
Then  \eqref{eq:etaIs}
explains the value of $\gd_{0}$ in \eqref{eq:gdAllow}.



Since the four regions cover the full range \eqref{eq:gagkRange}, we now just collect the results of the previous two sections.
In the range $\cR_{1}$, we apply Theorem \ref{thm:minBnd1},
 getting 
\be\label{eq:R1bnd}
\cI(Q,K)\ll {(\#\gW_{N})^{2}\over N} K^{-c}
\ee
for some $c>0$. In $\cR_{2}$, we apply Theorem \ref{thm:minBnd2}, 
getting
\be\label{eq:R2bnd}
\cI(Q,K)\ll {(\#\gW_{N})^{2}\over N} Q^{-c}
.
\ee
In the range $\cR_{3}$ with $K>\cQ$, Theorem \ref{thm:minBnd3} gives \eqref{eq:R1bnd}, and in $\cR_{4}$ with $Q>\cQ$, Theorem \ref{thm:minBnd4} gives \eqref{eq:R2bnd}.
Combining these estimates 
completes the proof of \eqref{eq:ErrBnd}. 
%
\epf

It is now standard to derive Theorem \ref{thm:mainN} from \eqref{eq:cMNbnd} and 
\eqref{eq:ErrBnd}.

\pf[Proof of Theorem \ref{thm:mainN}, assuming Proposition \ref{prop:cBII}]\

In light of \eqref{eq:cMNbnd}, the proof of which uses Proposition \ref{prop:cBII},  we have for $n\asymp N$ that 
\beann
\cM_{N}(n)
&\gg& 
|\gW_{N}|/(N\log\log N)
\\
&\gg&
 N^{2\gd-1-1/1000}
,
\eeann
where we crudely used \eqref{eq:gWLsize} and \eqref{eq:gbIs}.
Hence we expect the same for $R_{N}(n)$. If this is not the case, it means that 
$$
|\cE_{N}(n)|=|R_{N}(n)-\cM_{N}(n)|\gg {1\over\log\log N}{|\gW_{N}|\over N}
.
$$

Let $\fE(N)$ denote the set of $n\asymp N$  which have a 
small representation number
$R_{N}(n)$, 
$$
\fE(N):=
\left\{
\gba N\le n< \gbb N
:
R_{N}(n)< \foh\cM_{N}(n)
\right\}
.
$$
Then assuming \eqref{eq:gdBnd}, we have
\beann
\#\fE(N)
&\ll&
\sum_{
\gba N\le n< \gbb N
}\bo_{\{|\cE_{N}(n)|\gg{|\gW_{N}|\over N\log\log N}\}}
\\
&\ll&
{
N^{2}(\log\log N)^{2}
\over 
|\gW_{N}|^{2}}
\sum_{
n
}
|\cE_{N}(n)|
^{2}
\\
&\ll&
{
N^{2}(\log\log N)^{2}
\over 
|\gW_{N}|^{2}}
{
|\gW_{N}|^{2}\over N}
\cQ^{-c}
\ll
N^{1-c/\log\log N}
,
\eeann
using \eqref{eq:ErrBnd} and \eqref{eq:cQIs}.
\epf

This completes the proof of Theorem \ref{thm:mainN}, modulo the construction of the leading set $\aleph$, which is taken up in the next section.
First we give a quick

\subsection{Proof of Theorem \ref{thm:primes}}\label{sec:primes}

Let $\cP=\cP_{N}$ 
 be
  the set
 of primes $p$ up to $N$ which are $3(\mod4)$, 
so that
 $(p-1)/2$ 
is a $10$-almost-prime, that is,
 $$
\cP:=\{
p<N:p\equiv3(\mod 4)\text{, and }
m\mid(p-1)/2\Longrightarrow m>N^{1/10}
\}
.
 $$
A standard sieve argument shows that
$\cP$ 
  has 
cardinality $\gg
{N\over (\log N)^{2}}.
$
By \eqref{eq:smExcept},  the cardinality  $N^{1-c/\log\log N}$ of the exceptional set is much smaller, and so $\fD_{\cA}(N)\cap\cP$ is unbounded in $N$ for $\gd_{\cA}>\gd_{0}$.

By
\eqref{eq:mult}, each $p=d$ in the intersection  appears 
with multiplicity at least  $ N^{2\gd-1001/1000}> N^{10/11}$, say. 
That is, there are distinct $b_{1},\dots, b_{L}$ so that $b_{j}/d\in\fR_{\cA}$, $j=1,\dots,L$, and $L> N^{10/11}$.
Let $r$ be any primitive root mod $d$. For $j=1,\dots,L$, let $k_{j}$ be defined by $b_{j}\equiv r^{k_{j}}(\mod d)$, and let $K=\{k_{1},\dots,k_{L}\}$.
Of course $b_{j}$ is a primitive root mod $d$ iff $(k_{j},d-1)=1$.

Consider the subset $K'$ of $k\in K$ for which $(k,d-1)>2$. Since $d\in\cP$, each such $k$ has a prime factor of size $N^{1/10}$, and hence the cardinality of $K'$ is $\ll N^{9/10}$. This is less than the cardinality of $K$, so we may safely discard $K'$ from $K$, leaving a non-empty   set $K''$.

Consider $b\equiv r^{k}(\mod d)$ with $k\in K''$. If $(k,d-1)=1$, we are done, since $b$ is a root mod $d$ and $b/d\in\fR_{\cA}$. The only other possibility is $(k,d-1)=2$, whence $b$ is a square mod $d$.
Set $b':=d-b$, so $b'\equiv-r^{k}(\mod d)$; since $d\equiv 3(\mod 4)$, $b'$ is now a primitive root. It is elementary to verify that
 $b/d\in\fR_{\{1,2,\dots, A\}}$
 implies that
 $b'/d=1-b/d\in\fR_{\{1,2,\dots,A+1\}}$. That is, these quotients are still absolutely Diophantine, completing the proof. 

\section{Construction of  $\aleph$}\label{sec:New}

In this section, we
 arrange 
the special leading set $\aleph$ in the ensemble $\gW_{N}$ as described in \S\ref{sec:set0}.
We need two pieces of background, using \S\ref{sec:congTrans} to extract some modular/archimedean counting statements from \cite{BourgainGamburdSarnak2011}, and spending \S\ref{sec:randExt}
proving a certain ``randomness extraction argument.'' Finally, we proceed in \S\ref{sec:pfProp} to construct $\aleph$, thereby proving Proposition \ref{prop:cBII}, and finalizing the proof of  Theorem \ref{thm:mainN}.

\subsection{Congruence Counting Theorems%
}\label{sec:congTrans}


Recall from \S\ref{sec:Lalley} that
$\mu$
is
 the $\gd$-dimensional Hausdorff measure supported on the limit set $\fC$,
lifted to 
$\bP^{1}$. 
%
Extending the work of  Lalley \cite{Lalley1989}
to the congruence setting, 
Bourgain-Gamburd-Sarnak \cite{BourgainGamburdSarnak2011} proved the following theorem, adapted to our present context.

 \begin{thm}[\cite{BourgainGamburdSarnak2011}]
 \label{thm:BGS2}
There exists
an integer
\be\label{eq:fBis}
\fB=\fB(\cA)\ge1
\ee
and a constant
\be\label{eq:fcIs}
\fc=\fc(\cA)>0
\ee 
so that
the following holds.
For any $(q,\fB)=1$, any $\gw\in\SL_{2}(q)$, and any $\g_{0}\in\G$,
there is a constant $C(\g_{0})>0$ so that
\bea
\nonumber
&&
\hskip-1in
\#\left\{
\g\in\G
:
\g\equiv\gw(\mod q),
|v_{+}(\g)
-\fv|<\frac1H
,
\text{ and }
\frac{\|\g\g_{0}\|}{\|\g_{0}\|}\le T
\right\}
\\
\label{eq:thm31}
&&
\qquad
=
C(\g_{0})\cdot T^{2\gd}
{\mu(\cI)\over |\SL_{2}(q)|}
+
O\left(
T^{2\gd-\fc/\log\log T}
\right)
,
\quad
\text{ as $T\to\infty$.}
\eea
Here $\cI$ is the interval of length $1/H$ about $\fv$, and
the implied constant does not depend
 on $T$, $H$, $q$, $\gw$, or $\g_{0}$. 
Since $\fv$ is a density point and crudely using $\gd<1$, the main term in \eqref{eq:thm31} certainly exceeds the error if $H<q^{-3}T^{\fc/\log\log T}$.

With the same conditions as above, except for a modulus $q$ with $\fB\mid q$, we have
\bea
\label{eq:611s}
&&
\hskip-.5in
\#\left\{
\g\in\G
:
\g\equiv\gw(\mod q),
v_{+}(\g)\in\cI
,
\text{ and }
\frac{\|\g\g_{0}\|}{\|\g_{0}\|}\le T
\right\}
\\
\nonumber
&=&
{|\SL_{2}(\fB)|\over |\SL_{2}(q)|}
\cdot
\#\left\{
\g\in\G
:
\g\equiv\gw(\mod \fB),
v_{+}(\g)\in\cI
,
\text{ and }
\frac{\|\g\g_{0}\|}{\|\g_{0}\|}\le T
\right\}
\\
\nonumber
&&
\hskip3.5in
+
O(
T^{2\gd-\fc/\log\log T}
)
.
\eea
 \end{thm}

\begin{rmk}
Theorem \ref{thm:BGS2} was proved in \cite[see Theorem 1.5]{BourgainGamburdSarnak2011} under the further assumptions that the modulus $q$ is square-free, and that $\G$ is a convex-cocompact 
sub{\it group} of $\SL_{2}(\Z)$.
As discussed in \S\ref{sec:Lalley}, the proof is the same when the group is replaced by our free semigroup $\G$, which has no parabolic elements.
As for the level, taking $q$ square-free was enough for the sieving purposes in \cite{BourgainGamburdSarnak2011}, but for the circle method
used
 here, we must take arbitrary $q$. The main ingredient in analyzing the modular aspect (see \cite[Lemma 2]{BourgainGamburdSarnak2011}) was the spectral gap (expansion property) proved
using methods of 
additive combinatorics
and an $L^{2}$-flattening lemma
 in \cite{BourgainGamburd2008,BourgainGamburdSarnak2010}, again for square-free $q$. The relevant results have since been established in full generality for arbitrary modulus, see 
\cite{BourgainVarju2011} and \cite[Remark 30]{GolsefidyVarju2011}. With this input, \cite[Theorem 1.5]{BourgainGamburdSarnak2011} holds for arbitrary $q$.
With these two caveats, 
Theorem \ref{thm:BGS2} follows directly from the methods of \cite{BourgainGamburdSarnak2011}.
\end{rmk}

\begin{rmk}
Recall that throughout, the appearance of constants  $c$ and $C$ may change from line to line. 
The special constant $\fc$ in \eqref{eq:fcIs} is in contradistinction with this principle, being the same constant 
wherever it appears.
Moreover, Proposition \ref{prop:Lalley} follows immediately from \eqref{eq:thm31}, so the constant $\fc$ appearing there can be taken to be the same as the one here.
\end{rmk}

\begin{rmk}\label{rmk:fBmidQ}
In Theorem \ref{thm:BGS2}, we have stated the result only for the extreme cases $(\fB,q)=1$ and $\fB\mid q$. Of course intermediate cases can
be obtained by summing over suitable arithmetic progressions. This introduces no extra error since the number of terms is $\ll_{\fB}1$, and 
our implied constants may depend on the fixed alphabet $\cA$ (recall $\fB$ depends only on $\cA$).
\end{rmk}

The condition $\|\g\g_{0}\|/\|\g_{0}\|<T$ arises naturally in the above through the renewal method;
in fact, this condition is essentially equivalent to 
$$
d_{H}(\g\g_{0}i,i)-d_{H}(\g_{0}i,i)<C \log T
,
$$
where $d_{H}(z,w)$ denotes 
hyperbolic distance 
in the upper half plane $\bH$.
%
As in \S\ref{sec:Lalley}, one typically sets $\g_{0}=I$, but we will use $\g_{0}$ for a different purpose.
Namely, 
we will need to control both the expanding direction $v_{+}(\g)$, and its expanding eigenvalue $\gl(\g)$, but taking $\g_{0}=I$ gives us control only on
the norm $\|\g\|$ (which can be off by a constant from the eigenvalue). So we instead do the following.

It is easy to see from \eqref{eq:xie2e2} that
\be\label{eq:normToGl}
\|\g\|
=
{\gl(\g)
\over
|\<v_{+}(\g),v_{-}^{\perp}(\g)\>|}
\left(
1+
O\left(
\frac1{\|\g\|^{2}}
\right)
\right)
.
\ee
Assume now that both $\g$ and $\g_{0}$ lie in $\cI$, the interval of length $1/H$ about $\fv$, and assume $\|\g\|,\|\g_{0}\|>H$. Then applying \eqref{eq:normToGl}
to $\g_{0}$ and $\g\g_{0}$ gives
$$
\|\g_{0}\|
=
{\gl(\g_{0})
\over
|\<\fv,v_{-}^{\perp}(\g_{0})\>|}
\left(
1+
O\left(
\frac1{H}
\right)
\right)
,
$$
and
\beann
\|\g\g_{0}\|
&=&
{\gl(\g\g_{0})
\over
|\<v_{+}(\g\g_{0}),v_{-}^{\perp}(\g\g_{0})\>|}
\left(
1+
O\left(
\frac1{\|\g\g_{0}\|^{2}}
\right)
\right)
\\
&=&
{\gl(\g)\gl(\g_{0})
\over
|\<\fv,v_{-}^{\perp}(\g_{0})\>|}
\left(
1+
O\left(
\frac1{H}
\right)
\right)
,
\eeann
where we used \eqref{eq:1.7} and \eqref{eq:1.8}. On dividing, we obtain
$$
{
\|\g\g_{0}\|
\over
\|\g_{0}\|
}
=
{
\gl(\g)
}
\left(
1
+
O
\left(
\frac1H
\right)
\right)
,
$$
and can thus convert statements restricting norms into ones controlling  eigenvalues, without losing constants.

The constant $C(\g_{0})$ in \eqref{eq:thm31} approaches a constant $C(\fv)>0$ as $\|\g_{0}\|\to\infty$ with $v_{+}(\g_{0})\to\fv$; indeed, $C(\g_{0})$ is obtained by 
evaluating a certain Gibbs measure (see \cite[\S10]{BourgainGamburdSarnak2011} or \cite[(2.5)]{Lalley1989}, where his $x$ plays the role of our $\g_{0}$). Hence \eqref{eq:thm31} can be replaced by
\beann
&&
\hskip-1in
\#\left\{
\g\in\G
:
\g\equiv\gw(\mod q),
|v_{+}(\g)
-\fv|<\frac1H
,
\text{ and }
\gl(\g)\le T
\right\}
\\
&&
\qquad
=
C(\fv)\cdot T^{2\gd}
{\mu(\cI)\over |\SL_{2}(q)|}
\left(
1+
O\left(
\frac1H
\right)
\right)
+
O\left(
T^{2\gd-\fc/\log\log T}
\right)
,
\eeann
and a similar expression analogous to \eqref{eq:611s}.

Finally, we can restrict $\gl(\g)$ to a smaller range, $\gl(\g)=T(1+O(1/H_{1}))$, for a smaller parameter $H_{1}$; in applications we take $H_{1}=H^{1/2}$. In summary, we have the following.

\begin{cor}\label{cor:p7}
With notation as above, we have for any $T,H,H_{1}\to\infty$, and any
$(q, \fB)=1$, $\gw\in\SL_{2}(q)$ that
 \bea
\nonumber
&&
\hskip-1in
\#\left\{
\g\in\G
:
\g\equiv\gw(\mod q),
|v_{+}(\g)-\fv|<\frac1H
,
|\gl(\g)-T|<\frac T{H_{1}}
\right\}
\\
\label{eq:corp7}
&=&
C(\fv)\cdot {T^{2\gd} \over H_{1}}
{\mu(\cI)\over |\SL_{2}(q)|}
\left(
1+
O\left(
\frac1{H_{1}}
+
\frac{H_{1}}H
\right)
\right)
+
O\left(
T^{2\gd-\fc/\log\log T}
\right)
.
\eea
The implied constants are independent of $T$, $H$, $H_{1}$, $q$ and $\gw$.
If $H_{1}=o(H)$, then the main term dominates the error as long as $H_{1}H^{\gd+\vep}q^{3}\ll T^{\fc/\log\log T}$.

For a modulus $q\equiv0(\mod \fB)$, we have
\bea
\label{eq:cor2}
&&
\hskip-.5in
\#\left\{
\g\in\G
:
\g\equiv\gw(\mod q),
v_{+}(\g)\in\cI
,
|\gl(\g)-T|<\frac T{H'}
\right\}
\\
\nonumber
&=&
{|\SL_{2}(\fB)|\over |\SL_{2}(q)|}
\cdot
\#\left\{
\g\in\G
:
\g\equiv\gw(\mod \fB),
v_{+}(\g)\in\cI
,
|\gl(\g)-T|<\frac T{H'}
\right\}
\\
\nonumber
&&
\hskip1in
\times
\left(
1+
O\left(
\frac1{H_{1}}
+
\frac{H_{1}}H
\right)
\right)
\quad+\quad
O(
T^{2\gd-\fc/\log\log T}
)
.
\eea
\end{cor}

%

\subsection{A Randomness Extraction Argument}\label{sec:randExt}\

Corollary \ref{cor:p7} 
gives us good modular/achimedean control away from the modulus $\fB$, but we need $\aleph$ to have good distribution properties for all moduli. So we will concoct in the next subsection certain special sets engineered to have good equidistribution mod $\fB$. But in so doing, we will potentially ruin the distribution away from $\fB$. To recover this distribution, we apply a certain
more-or-less standard
 ``randomness extraction'' argument, which states roughly that 
 if a large set has good modular distribution,
then so does a sufficiently large random subset of it. 
We will need to have the flexibility to stay away from a modulus $q_{0}$, which in applications is either $1$ or $\fB$.

\begin{lem}\label{lem:ranEx}
Let $\mu=\mu_{S}$ be the normalized (probability) measure of a finite subset 
$S\subset\SL(2,\Z)$, 
$$
\mu(\g)=\frac1{|S|}\sum_{s\in S}\bo_{\{s=\g\}}
,
$$
and fix $\eta>0$.
Let $q_{0}<Q$ be
a fixed modulus,
 let $\gw_{0}\in\SL_{2}(q_{0})$ be a fixed element,
and
let $\fQ=\fQ_{q_{0}}\subset[1,Q]$ be the set of moduli $q<Q$ with $q_{0}\mid q$. 
Assume that for all
$q\in\fQ$ and all
 $\gw\in\SL_{2}(q)$ with $\gw\equiv\gw_{0}(\mod q_{0})$, the projection 
$$
\pi_{q}[\mu](\gw)
=
\sum_{\g\equiv\gw(\mod q)}\mu(\g)
$$ 
is near the uniform measure on $\SL_{2}(q)$ conditioned on being $\equiv\gw_{0}(\mod q_{0})$,
\be\label{eq:piqMu}
\left\|
\pi_{q}[\mu]-\frac{|\SL_{2}(q_{0})|}{|\SL_{2}(q)|}
\right\|
_{L^{\infty}\big|_{\equiv\gw_{0}(\mod q_{0})}}
=
\max_{\gw\in\SL(2,q)\atop \gw\equiv\gw_{0}(\mod q_{0})}
\left|
\pi_{q}[\mu](\gw)-\frac{|\SL_{2}(q_{0})|}{|\SL_{2}(q)|}
\right|
<
\eta.
\ee
Then for any 
\be\label{eq:Tsize}
T>\eta^{-2}\log Q
,
\ee 
there exist $T$ distinct points $\g_{1},\dots, \g_{T}\in S=\supp\mu$ 
such that the probability measure  $\nu=\nu_{T,\g_{1},\dots,\g_{T}}$ defined by
\be\label{eq:nuDef}
\nu=\frac1T\left(\bo_{\g_{1}}+\cdots+\bo_{\g_{T}}\right)
\ee
has 
the same property.
That is, for all $q\in\fQ$ projection $\pi_{q}[\nu]$ is also nearly uniform,
\be\label{eq:piqNu}
\max_{q\in\fQ}
\left(
\left\|\pi_{q}[\nu]-\frac{|\SL_{2}(q_{0})|}{|\SL_{2}(q)|}\right\|
_{L^{\infty}\big|_{\equiv\gw_{0}(\mod q_{0})}}
\right)
\ll\eta.
\ee
The  implied constant above is absolute.
\end{lem}

\pf
This is a standard argument, so we give a sketch. 
Take $\nu$ as in \eqref{eq:nuDef}.
Let $\cD$ be the 
expectation with respect to $\mu$ of the 
left hand side of \eqref{eq:piqNu},
\beann
\cD
&=&
\max_{q\in\fQ}
\sum_{\g\in\otimes^{T}\SL_{2}(\Z)}
\max_{\gw\in\SL_{2}(q)\atop\gw\equiv\gw_{0}(\mod q_{0})}
\left|
\frac1T
\sum_{j=1}^{T}
\bo_{\{\g_{j}\equiv\gw(q)\}}
-
\frac{|\SL_{2}(q_{0})|}{|\SL_{2}(q)|}
\right|
\mu^{(T)}(\g)
,
\eeann
where $\mu^{(T)}$ is the product measure on $\otimes^{T}\SL_{2}(\Z)$ and $\g=(\g_{1},\dots,\g_{T})$.
Using \eqref{eq:piqMu}, we have
\beann
\cD
&<&
\eta
+
\max_{q\in\fQ}
\sum_{\g\in\otimes^{T}\SL_{2}(\Z)}
\sum_{\xi\in\otimes^{T}\SL_{2}(\Z)}
\max_{\gw\in\SL_{2}(q)
}
\left|
\frac1T
\sum_{j=1}^{T}
f_{\gw}(\g_{j},\xi_{j})
\right|
\mu^{(T)}(\g)
\mu^{(T)}(\xi)
,
\eeann
where
$$
f_{\gw}(\g_{j},\xi_{j})
:=
\bo_{\{\g_{j}\equiv\gw(q)\}}-\bo_{\{\xi_{j}\equiv\gw(q)\}}
,
$$
and we extended the max over $\gw$ to all of $\SL_{2}(q)$.

Note that for fixed $\gw$, $f_{\gw}(\g_{j},\xi_{j})$ are independent, mean zero random variables and bounded by $1$. Hence the contraction principle gives
\be\label{eq:cDbackTo}
\cD
<
\eta
+
\max_{q\in\fQ}
\sum_\g
\sum_\xi
\cD_{q}(\g,\xi)
\mu^{(T)}(\g)
\mu^{(T)}(\xi)
,
\ee
where
$$
\cD_{q}(\g,\xi)
:=
\frac1{2^{T}}
\sum_{\vep\in\{\pm1\}^{T}}
\max_{\gw\in\SL_{2}(q)
}
\left|
\frac1T
\sum_{j=1}^{T}
\vep_{j}
f_{\gw}(\g_{j},\xi_{j})
\right|
.
$$
Replace the $\max$ by an $L^{p}$ norm with $p$ to be chosen later:
\beann
\cD_{q}(\g,\xi)
&\le&
\frac1{2^{T}}
\sum_{\vep\in\{\pm1\}^{T}}
\left(
\sum_{\gw\in\SL_{2}(q)}
\left|
\frac1T
\sum_{j=1}^{T}
\vep_{j}
f_{\gw}(\g_{j},\xi_{j})
\right|^{p}
\right)^{1/p}
\\
&\ll&
\left(
\sum_{\gw\in\SL_{2}(q)}
\frac1{2^{T}}
\sum_{\vep\in\{\pm1\}^{T}}
\left|
\frac1T
\sum_{j=1}^{T}
\vep_{j}
f_{\gw}(\g_{j},\xi_{j})
\right|^{p}
\right)^{1/p}
\\
&\ll&
\left(
\sum_{\gw\in\SL_{2}(q)}
p^{p/2}
\left(
\sum_{j=1}^{T}
\left|
\frac
{f_{\gw}(\g_{j},\xi_{j})}
T
\right|^{2}
\right)^{p/2}
\right)^{1/p}
\\
&\ll&
q^{3/p}
p^{1/2}
T^{-1/2}
,
\eeann
where we applied 
Khintchine's 
inequality \cite{Haagerup1981}
(the implied constant is absolute). 
Now we choose $p=\log q$, so that
$$
\cD_{q}
\ll
(\log q)^{1/2}
T^{-1/2}
\le
(\log Q)^{1/2}
T^{-1/2}
.
$$
Inserting this into \eqref{eq:cDbackTo} and setting $T>\eta^{-2}\log Q$ gives
$$
\cD\ll \eta,
$$
from which the claim follows immediately.
\epf

Equipped with this randomness extraction argument, we proceed with the 
\subsection{Proof of Proposition \ref{prop:cBII}}\label{sec:pfProp}
\

Recalling the parameters $\Bpow$ in  \eqref{eq:gbIs}, $\fc$ from \eqref{eq:fcIs}, and $\fB$ from \eqref{eq:fBis}, we  set $R:=|\SL_{2}(\fB)|$, and 
 define $\ga_{0}$ in \eqref{eq:cQIs} by
\be\label{eq:ga0Is}
\ga_{0}:={\gb\fc\over 40 R}
.
\ee 
For a parameter 
\be\label{eq:TNc}
T=N^{c_{1}}
\ee 
with small $c_{1}$ to be determined in \eqref{eq:TruDef}, let $H=\cQ^{12}$, $H_{1}=\cQ^{6}$, and set
\be\label{eq:cSTdef}
\cS(T):=
\left\{
\g\in\G
:
|v_{+}(\g)-v|<\frac1H
,
|\gl(\g)-T|<\frac T{H_{1}}
\right\}
.
\ee
By \eqref{eq:corp7} with $q=1$, we have
crudely using $\gd<1$ that
\be\label{eq:STcount}
\#\cS(T)
\gg
T^{2\gd}/\cQ^{18}
+
O(T^{2\gd-\fc/\log\log T})
.
\ee
We have from  \eqref{eq:TNc}  that 
$$
T^{-\fc/\log\log T} 
\ll 
T^{-\fc/\log\log N} 
\ll 
\left(
N^{\fc/\log\log N} 
\right)^{-c_{1}}
\ll 
(\cQ^{40R/\gb})^{-c_{1}}
.
$$
So as long as 
\be\label{eq:TpowIs}
c_{1}>\TwenOvQpow
,
\ee
\eqref{eq:STcount} is significant, with an error of size $\ll T^{2\gd}/\cQ^{30}$. 

By the pigeonhole principle, there is some element $\fs_{T}\in\cS(T)$ so that
$$
\cS'(T):=
\{
s\in \cS(T)
:
s\equiv \fs_{T}(\mod \fB)
\}
$$
satisfies
\be\label{eq:cSpTcount}
\#\cS'(T)
\ge
{1\over |\SL_{2}(\fB)|}
\#\cS(T)
\gg
T^{2\gd}/\cQ^{18}
.
\ee
(Recall our implied constants may depend implicitly on $\cA$, and $\fB$ depends only on $\cA$.)

For this set, the counting statement \eqref{eq:STcount} remains significant even with a modular restriction: for any $q<\cQ$ with $\fB\mid q$, and any $\gw\in\SL_{2}(q)$ with $\gw\equiv\fs_{T}(\mod \fB)$, applying \eqref{eq:cor2} gives 
\bea
\nonumber
&&
\hskip-.5in
\#
\{
s\in\cS'(T)
:
s\equiv \gw(\mod q)
\}
=
\#
\{
s\in\cS(T)
:
s\equiv \gw(\mod q)
\}
\\
\nonumber
&=&
{|\SL_{2}(\fB)|\over |\SL_{2}(q)|}
\#
\{
s\in\cS(T)
:
s\equiv \gw\equiv\fs_{T}(\mod \fB)
\}\
(1+O(\cQ^{-6}))
+O(T^{2\gd}\cQ^{-30})
\\
\label{eq:SpModq}
&=&
{|\SL_{2}(\fB)|\over |\SL_{2}(q)|}
\#
\cS'(T)
\ (1+O(\cQ^{-6}))
+
O(T^{2\gd}\cQ^{-30})
.
\eea
From \eqref{eq:cSpTcount}, the main term is $\gg T^{2\gd}/\cQ^{21}$, dominating the error.

We just need to play with $\cS'(T)$ to get good distribution modulo $\fB$. 
Recall $R=|\SL_{2}(\fB)|$. 
Then 
every element of the ``coset'' 
$$
\g\in\cS'(T)\cdot \fs_{T}^{R-1}
$$ 
satisfies $\g\equiv I(\mod \fB)$.
Next write
$
\SL_{2}(\fB)
=
\{\g_{1},\g_{2},\dots,\g_{R}\}$, and take 
$x_{1},\dots, x_{R}\in\G$ so that
\be\label{eq:3.14}
x_{r}\equiv \g_{r}(\mod\fB)
,\qquad
r=1,\dots,R.
\ee
(Recall we had assumed in \S\ref{sec:set0}  that $\G(\mod q)$ is all of $\SL_{2}(q)$ for all $q$, so such $x_{r}$ exist.)
Such $x_{r}$ can be found of size $\asymp_{\cA}1$. 

Note that any  element 
$$
\g\in\cS'(T)\cdot \fs_{T}^{R-1}\cdot x_{r}
$$ 
has $\g\equiv \g_{r}(\mod \fB)$. 
Unfortunately, this triple-product does not work, since 
we do not have control on the expanding
vector of $x_{r}$. 
To remedy this,
we take a single fixed element $\ff_{0}\in\G$ of size
\be\label{eq:ffsize}
\gl(\ff_{0})\asymp B^{1/100},
\ee 
say, 
with 
\be\label{eq:ffvp}
|v_{+}(\ff_{0})-\fv|<\cQ^{-6}
.
\ee 
Then
from \eqref{eq:1.8},
 $v_{+}(\ff_{0}x_{r})=\fv(1+O(\cQ^{-6}))$, and for any $s\in\cS'(T)$, 
$$
v_{+}(s\cdot\fs_{T}^{R-1}\cdot\ff_{0}x_{r})=\fv(1+O(\cQ^{-6}))
.
$$
Moreover from \eqref{eq:cSTdef} and \eqref{eq:1.7}, we have
\beann
\gl(s\cdot\fs_{T}^{R-1}\cdot\ff_{0}x_{r})
&=&
\gl(s)\gl(\fs_{T})^{R-1}\gl(\ff_{0}x_{r})
\left(1+O(\cQ^{-6})\right)
\\
&=&
T^{R}\gl(\ff_{0}x_{r})
\left(1+O(\cQ^{-6})\right)
.
\eeann

Now for each $u\in\cU$, $u\asymp B$, and each $r=1,\dots,R$, take $T=T_{u,r}$ so that
$$
T^{R}\gl(\ff_{0}x_{r})
=u
,
$$
that is, let
\be\label{eq:TruDef}
T_{u,r}
:=
\left({u\over\gl(\ff_{0}x_{r})}\right)^{1/R}
\asymp
B^{99/(100R)}
=
N^{99\Bpow/(100R)}
,
\ee
which, by \eqref{eq:Bis}, determines $c_{1}$ in \eqref{eq:TNc}.
Note that \eqref{eq:TpowIs} is easily satisfied. 

Thus for each $u$ and $r$, we have sets
$$
\cB_{u,r}:=\cS'(T_{u,r})\cdot (\fs_{T_{u,r}})^{R-1}\cdot \ff_{0}\cdot x_{r}\subset\G
,
$$
so that for all $\fa\in\cB_{u,r}$, the expanding vector is controlled,
$$
|v_{+}(\fa)-v|\ll \cQ^{-6},
$$ 
and the eigenvalue is controlled,
$$
\gl(\fa)=u(1+O(\cQ^{-6}))
.
$$
Since we have saved an extra $\cQ$, we can use it to set the implied constant to 1, getting \eqref{eq:d8p} and \eqref{eq:d9p}.

 Note
 that
 by  \eqref{eq:SpModq},
  for all $q<\cQ$ with $\fB\mid q$, and all $\gw\in\SL_{2}(q)$ 
with $\gw\equiv\ff_{0}x_{r}(\mod \fB)$, 
we have, crudely, that
\be\label{eq:BurModq}
\#\{
\fa\in\cB_{u,r}:
\fa\equiv\gw(\mod q)
\}
=
{|\SL_{2}(\fB)|\over|\SL_{2}(q)|}
\#
\cB_{u,r}
\left(
1+O(\cQ^{-5})
\right)
.
\ee

Recall also from \eqref{eq:STcount} that the cardinality of $\cB_{u,r}$ is 
\be\label{eq:BurSize}
\gg (T_{u,r})^{2\gd}/\cQ^{18}
\gg N^{c}
,
\ee
using \eqref{eq:TruDef}.

Hence for fixed $u$, we may apply the randomness extraction argument in Lemma \ref{lem:ranEx} 
to $\cB_{u,r}$,  with $\eta=\cQ^{-5}$ and
$q_{0}=\fB$.
%
%
This gives sets $\cB'_{u,r}\subset\cB_{u,r}$ of size $\gg N^{c}$, for which 
\eqref{eq:BurModq} continues to hold, and moreover we can force them all to have exactly the same cardinality independently of $r$,
$$
|\cB'_{u,r}|
=
|\cB'_{u,r'}|
.
$$
Set
$$
\tilde\aleph_{u}:=\bigsqcup_{r=1}^{R}\cB'_{u,r},
$$
and note that for $\fB\mid q<\cQ$ and $\gw\in\SL_{2}(q)$,
\beann
\#\{
\fa\in
\tilde\aleph_{u}
:
\fa\equiv\gw(\mod q)
\}
=
{|\SL_{2}(\fB)|\over |\SL_{2}(q)|}
\#
\cB'_{u,r}
\left(1+O(\cQ^{-5})\right)
,
\eeann
where $r$ is the index for which $\gw\equiv \ff_{0}x_{r}(\mod \fB)$.
Since
$
\#
\cB'_{u,r}
=
{|\tilde\aleph_{u}|/R}
,
$
and
$R=|\SL_{2}(\fB)|$,
we have that
for each $u$, $\tilde\aleph_{u}$ satisfies 
\be\label{eq:tilaleph}
\#\{
\fa\in
\tilde\aleph_{u}
:
\fa\equiv\gw(\mod q)
\}
=
{|\tilde\aleph_{u}|\over |\SL_{2}(q)|}
\left(1+O(\cQ^{-5})\right)
.
\ee
We can now also drop the condition $\fB\mid q$ in \eqref{eq:tilaleph} by summing along certain arithmetic progressions; since $\fB\ll_{A}1$, the implied constant still depends only on $A$, cf. Remark \ref{rmk:fBmidQ}.

Now we apply Lemma \ref{lem:ranEx} again to $\tilde\aleph_{u}$, with $\eta=\cQ^{-5}$ and $q_{0}=1$, giving sets 
$$
\aleph_{u}\subset\tilde\aleph_{u}
$$ 
for which 
\eqref{eq:tilaleph}
still holds, that is \eqref{eq:cBSatp} holds, and which all have the same cardinality, giving \eqref{eq:alephUsize}.

This completes the proof of Proposition \ref{prop:cBII}.



\section{Proof of Theorem \ref{thm:gdPlus}}\label{sec:appendix}

Recall that $\fR_{\cA}(N)$ is the set of rationals $b/d$ with partial quotients bounded in the alphabet $\cA$ with $d<N$,  $\fD_{\cA}(N)$ is the set of continuants up to $N$, and $\fC_{\cA}$ is the limit set of $\fR_{\cA}$ with Hausdorff dimension $\gd=\gd_{\cA}$. Recall the sum-set structure \eqref{eq:sumSet}, that if $a\in\cA$ and $b/d\in\fR_{\cA}$, then $d$ and $b+ad$ are in $\fD_{\cA}$. The same holds for another $a'\in\cA$, that is, all three of $d$, $b+ad$, and  $b+a'd$ are in $\fD_{\cA}$.

We wish to show \eqref{eq:gdPlus} that for any $\vep>0$,
$$
\#(\fD_{\cA}\cap[1,N])\gg_{\vep}N^{\gd+(2\gd-1)(1-\gd)/(5-\gd)-\vep}
.
$$
For ease of notation, we lose no generality by specializing from now on to the case $1,2\in\cA$; whence $b/d\in\fR_{\cA}(N)$ implies that
\be\label{eq:ap1.6}
d,b+d,b+2d\in\fD_{\cA}(3N)
.
\ee
And again, we can drop the  subscript $\cA$ from 
$\fD$, $\fR$ and $\fC$. 

Our new ensemble of focus is the collection $\gW=\gW_{N}$ of intervals given by
\be\label{eq:ap2.2}
\gW:=
\bigcup_{d\in\fD\atop N/2\le d<3N}\left[\frac {d}N,\frac{d+1}N\right]
\subset[1/2,3]
,
\ee
so that
\be\label{eq:gWtofD}
|\gW|\ll {\#(\fD\cap[1,3N])\over N}
.
\ee
Hensley's conjecture implies, assuming $\gd>1/2$, that we should have 
\be\label{eq:gW1}
|\gW|\asymp1.
\ee
A priori, we do not even know that $|\gW|>1$, and the bound \eqref{eq:DANpower} follows from
\be\label{eq:gWN1md}
|\gW|\gg N^{-(1-\gd)}
,
\ee
so this is what we must beat.

Note from \eqref{eq:ap1.6} that 
\be\label{eq:bogW}
\bo_{\gW}(x)\bo_{\gW}(y)\bo_{\gW}(x+y)\bo_{\fR}\left(\frac yx-1\right)
=1
\ee
if $x=d/N$ and $y=(b+d)/N$ with $b/d\in\fR$. 
Just as we thickened $\bo_{\fD}$ to some intervals $\bo_{\gW}$ in \eqref{eq:ap2.2}, we wish  to thicken $\bo_{\fR}$ 
to some intervals.

By 
Frostman's
theorem,
there is a probability measure
 $\mu$
supported  on the Cantor set $\fC$, 
so that
for any interval $\cI\subset[0,1]$, we have
\be\label{eq:ap1.1}
\mu(\cI)\ll_{\vep} |\cI|^{\gd-\vep}
.
\ee
The most naive thickening of $\fR$ we could take is at the scale of $1/N^{2}$; namely, for each $x\in\fC$, there is (by Dirichlet's approximation theorem and properties of continued fractions) some $b/d\in\fR(N)$ with 
$$
\left|x-\frac bd\right|<\frac1{N^{2}}
.
$$ 
So we can find a collection of $\asymp N^{2\gd}$ intervals of length $\asymp 1/N^{2}$ which cover $\fC$, each of which has a point in $\fR(N)$. Note also that the spacing between consecutive points in $\fR(N)$ is $\ge 1/N^{2}$. Instead it will be more fruitful to collect points at square-root this scale, $1/N$, as follows.

By \eqref{eq:ap1.1}, there is 
a collection
 $\{\cI_{\ell}\}_{\ell\le L}$
of 
$L\asymp N^{\gd}$
 disjoint intervals $\cI_{\ell}\subset[1/2,1]$, each of length $1/N$,
 so that 
\be\label{eq:ap1.2}
\mu(\cI_{\ell})\gg_{\vep} N^{-\gd-\vep}
.
\ee
Denote their union by
\be\label{eq:tildeR}
\tilde\cR:=\bigsqcup_{\ell=1}^{L}\cI_{\ell}
.
\ee
Subdividing each $\cI_{\ell}$ into intervals $\cI_{\ell,n}$ of length $1/N^{2}$, it follows that there are at least $\gg N^{\gd-\vep}$ of them intersecting $\fC$.
For each such $\ell,n$, 
the intersection
$
\cI_{\ell,n}\cap\fR(N)$ is also non-empty, possibly after replacing $\cI_{\ell,n}$ with a doubling. 
Hence
the cardinality of
$$
\tilde\cR(N):=
\left\{\frac bd\in\fR\cap\tilde\cR:(b,d)=1,
N/2<b<d<N\right\}
$$
is of the right order
\be\label{eq:tilRNbnd}
\#\tilde\cR(N)
\gg_{\vep} N^{2\gd-\vep}
.
\ee 
We thicken these intervals just a little further, setting
\be\label{eq:ap2.1}
\cR:=\bigcup_{\ell=1}^{L}\left(\cI_{\ell}+\left[-\frac2N,\frac2N\right]\right)
.
\ee
Note that since $L\ll N^{\gd}$, we have
\be\label{eq:cRsize}
|\cR|\ll N^{-(1-\gd)}
.
\ee
Note further that $\fC$, and hence $\cR$, is contained strictly inside $[0,1]$, that is, for some $\nu=\nu_{A}>0$, we have
\be\label{eq:apRcont}
\cR\subset[\nu,1-\nu].
\ee

Consider now our main integral $\cJ$, motivated by \eqref{eq:bogW}, defined by
\be\label{eq:cJdef}
\cJ:=
\iint_{\R^{2}}
\bo_{\gW}(x)
\bo_{\gW}(y)
\bo_{\gW}(x+y)
\bo_{\cR}\left(\frac yx-1\right)
dydx
.
\ee
By \eqref{eq:ap2.2}, the domain of integration above is supported in the box $[1/2,3]\times[1/2,3]$.

For each $b/d\in\tilde\cR(N)$,  the intervals
$$
{d\over N}
\le x \le
{d+1/2\over N}
,
\qquad
{b+d\over N}
\le y\le {b+d+1/2\over N}
,
$$ 
belong to $\gW$, by \eqref{eq:ap1.6}, as does the interval
$$
{b+2d\over N}
\le x+y \le
{b+2d+1\over N}
.
$$
Moreover the interval
$$
\qquad
{b\over d}-\frac2N<{b-1\over d}
\le \frac yx-1\le {b+1\over d}
<\frac bd+\frac2N
$$
is in $\cR$ by the thickening in \eqref{eq:ap2.1}.
That is, these intervals contribute $1/N^{2}$ to $\cJ$ for each $b/d$, and hence
by \eqref{eq:tilRNbnd}, we have established the following lower bound.
\begin{prop}
\be\label{eq:ap2.4}
\cJ\gg_{\vep} 
N^{-2(1-\gd)-\vep}
.
\ee
\end{prop}

The trivial bound $\cJ\ll|\gW|^{2}$ from \eqref{eq:cJdef}
recovers \eqref{eq:gWN1md
}, which we want to beat.
The rest of the appendix is devoted to establishing the following

\begin{prop}\label{prop:Jupper}
\be\label{eq:Jupper}
\cJ\ll_{\vep} N^{-{1-\gd\over2}+\vep}|\gW|^{2- {3(1-\gd)\over2(2-\gd)}}
.
\ee
\end{prop}

Then Theorem \ref{thm:gdPlus} follows immediately from \eqref{eq:ap2.4}, \eqref{eq:Jupper}, and \eqref{eq:gWtofD}.

\subsection{Proof of Proposition \ref{prop:Jupper}}\

Let $M>0$ be a parameter to be chosen later (it will be a little less than $N^{-(1-\gd)}$), and decompose
\be\label{eq:cJsplit}
\cJ=\cJ_{1}+\cJ_{2}
,
\ee
where
\be\label{eq:ap3.1}
\cJ_{1}:=\iint_{(\bo_{\gW}\ast \bo_{-\gW})(-x)<M}\cdots dxdy
,
\ee
and
\be\label{eq:ap3.2}
\cJ_{2}:=\iint_{(\bo_{\gW}\ast \bo_{-\gW})(-x)\ge M }\cdots dxdy
.
\ee

Then writing $\bo_{\cR}\le1$, we have 
\bea
\nonumber
\cJ_{1}
&\le&
\iint\limits_{(\bo_{\gW}\ast \bo_{-\gW})(-x)<M}
\bo_{\gW}(x)
\bo_{\gW}(y)
\bo_{-\gW}(-x-y)
dxdy
\\
\nonumber
&\le&
\int\limits_{(\bo_{\gW}\ast \bo_{-\gW})(-x)<M}
\bo_{\gW}(x)
\bigg(
\bo_{\gW}\ast 
\bo_{-\gW}(-x)
\bigg)
dx
\\
\label{eq:ap3.3}
&<&
M
|\gW|
.
\eea
It is clear already that to get a gain on $|\gW|$, we must take $M$ a power less than $N^{-(1-\gd)}$. \\

We are left to analyze $\cJ_{2}$. 
Note that in the  domain of $\cJ_{2}$, we have
$$
\bo_{\gW}(x+y)
\le1\le 
\frac1M(\bo_{\gW}\ast \bo_{-\gW})(-x)
.
$$ 
Hence we can write
$$
\cJ_{2}
\le
\frac1M
\iint_{\R^{2}}
\eta(x)\
(\bo_{\gW}\ast \bo_{-\gW})(-x)\
\bo_{\gW}(y)\
\bo_{\cR}\left(\frac yx-1\right)
dy
dx
,
$$
where we
have bounded $\bo_{\gW}(x)$ by
a 
smooth 
bump function $\eta(x)$ 
 with support in $[1/4,4]$, say, and $\eta\ge1$ on $[1/2,3]$ to recall \eqref{eq:ap2.2} that $x\in[1/2,3]$.

For a smooth, non-negative, even function $\gU$ with compact support and $\int \gU=1$, let $\gU_{N}(y):=10N\gU(10Ny)$, and dominate $\bo_{\gW}(y)$, up to constant, by the smooth function
$$
\cS_{\gW}:=\bo_{\gW}\ast \gU_{N}
.
$$
So we have
\be\label{eq:J2J2p}
\cJ_{2}
\ll
\frac1M
\cdot \cJ'_{2},
\ee
where
\be\label{eq:ap3.4}
\cJ'_{2}
:=
\iint_{\R^{2}}
\eta(x)
(\bo_{\gW}\ast \bo_{-\gW})(-x)\
\cS_{\gW}(y)\
\bo_{\cR}\left(\frac yx-1\right)
dy
dx
.
\ee

\begin{figure}
%
\includegraphics[width=5in]{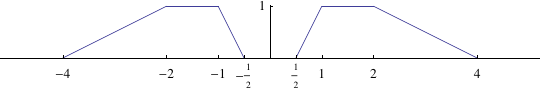}
%
%
%
\caption{
The Fourier multiplier $\gl_{0}$.
}
\label{fig:gl}
\end{figure}

Note that, by the smoothness of $\gU$, the Fourier spectrum of $\cS_{\gW}$ is contained, up to negligible error, in $[-N^{1+\vep},N^{1+\vep}]$. 
So  we can decompose $\cS_{\gW}(y)$ by the technique of ``slicing.'' That is, introduce a certain dyadic partition of unity via the 
Fourier multipliers $\gl_{k}(\xi)$,
defined
as follows. 
Let $\gl_{0}(\xi)$ be even, $\equiv1$ on $[1,2]$, and decaying piecewise-linearly to $0$ at $\xi=1/2$ and $\xi=4$,
see Figure \ref{fig:gl}. 
For integers $k$ ranging in
\be\label{eq:apkBnd}
0<2^{2k}<N^{1+\vep}
,
\ee
define
$$
\gl_{k}(\xi):=\gl_{0}(\xi\cdot2^{-2k}).
$$
%
%
%
Let $\gL_{k}$ be the Fourier inverse of $\gl_{k}$, so 
$$
\gL_{k}(y)
=2^{2k}\gL_{0}(2^{2k}y)
=(\cD_{2^{2k}}\gL_{0})(y)
,
$$ 
where
$$
\gL_{0}(y):=
\frac{ \sin ^2\left(\frac32 \pi   y\right) }{\pi ^2 y^2}
\bigg(
   2 \cos \left(2\pi  
   y\right)
-\cos \left(\pi  
   y\right)
   +\cos \left(5 \pi   y\right)
   \bigg)
,
$$
and 
 $\cD_{u}$ is the dilation representation,
$$
(\cD_{u}f)(y):=uf(uy)
.
$$
We also introduce  $\cT_{u}$, the translation representation, 
$$
\cT_{u}f(y):=f(y+u)
.
$$

Of course $\widehat{\gL_{k}\ast f}=\gl_{k}\cdot\hat f$, so we have
$$
\cS_{\gW}=\sum_{
2^{2k}<N^{1+\vep}
}(\gL_{k}\ast \cS_{\gW})
+
Err
,
$$
where $Err$ is bounded by an arbitrarily large power of $1/N$, and will henceforth be ignored.

Then we can bound \eqref{eq:ap3.4} as
\be\label{eq:cJ2k}
\left|
\cJ'_{2}
\right|
\ll
\sum_{k
}
\left|
\cJ_{2}^{(k)}
\right|
,
\ee
where
\be\label{eq:cJkIs}
\cJ_{2}^{(k)}
:=
\iint
\eta(x)
\
(\bo_{\gW}\ast \bo_{-\gW})(-x)\
(\gL_{k}\ast\cS_{\gW})(y)\
\bo_{\cR}\left(\frac yx-1\right)
dy
dx
.
\ee

Since the Fourier spectrum of $\cS_{\gW}$ is now controlled, so is that of $\bo_{\cR}$, as follows.
Write 
$
\bo_{\cR}\left(\frac yx-1\right)
=
\bo_{x\cR}\left( y-x\right)
=
(\cT_{-x}\bo_{x\cR})(y)
.
$
Then the $y$ integral can be written as
\bea
\label{eq:glkp}
&&
\hskip-.4in
\int_{\R}
(\gL_{k}\ast\cS_{\gW})(y)\
(
\cT_{-x}\bo_{x\cR})
(y)
dy
\\
\nonumber
&&
=
\int_{\R}
\gl_{k}(\xi)\hat\cS_{\gW}(\xi)\
\overline{
\widehat{
(
\cT_{-x}\bo_{x\cR})
}(\xi)}
d\xi
\\
\nonumber
&&
=
\int_{\R}
\gl_{k}(\xi)\hat\cS_{\gW}(\xi)\
\gl'_{k}(x\xi)
\overline{
\widehat{
(
\cT_{-x}\bo_{x\cR})
}(\xi)}
d\xi
,
\eea
where we inserted
another bump function
 $\gl'_{k}$
which is smooth in addition to other properties of $\gl_{k}$. Namely, let $\gl_{0}'$ be even, $\equiv1$ on $\pm[1/16,16]$, and decay smoothly to $0$ outside of $\pm[1/32,32]$; then set $\gl'_{k}(\xi):=\gl'_{0}(\xi 2^{-2k})$. The point is that $\gl'_{k}(x\xi)\equiv1$ on the support of $\gl_{k}$, since $x\in[1/4,4]$ by the support of $\eta$, so the above equality holds.

 Then writing $\gL'_{k}$ for the inverse transform of $\gl'_{k}$,
  we have 
\be
\label{eq:ap3.5}
\cJ_{2}^{(k)}
=
\iint
\eta(x)\
(\bo_{\gW}\ast \bo_{-\gW})(-x)\
(\gL_{k}\ast\cS_{\gW})(y)\
\bigg[
\gL'_{k}
\ast
\bo_{\cR}
\bigg]
\left(\frac yx-1\right)\
dy\,
dx
.
\ee

Now we handle two ranges of $k$ separately. We introduce a cutoff parameter $\cK
$ to be chosen later, see \eqref{eq:cKIs}.

\subsubsection{The range $k\le \cK$}\

We wish to prove
\begin{lem}
For $k\le \cK$,
\be\label{eq:ap3.8}
\left|\cJ_{2}^{(k)}\right|
\ll_{\vep}
|\gW|^{3}\
2^{2k(1-\gd)}\
N^{-(1-\gd)+\vep}
.
\ee
\end{lem}

This is a gain of a power of $|\gW|$ (recall we are assuming $|\gW|<1$).

\pf
Estimate \eqref{eq:ap3.5} by
\be\label{eq:J2k1}
\left|\cJ_{2}^{(k)}\right|
\ll
\left\|\bo_{\gW}\ast \bo_{-\gW}\right\|_{1}
\cdot
\left\|\gL_{k}\ast\cS_{\gW}\right\|_{1}
\cdot
\bigg\|
\gL'_{k}\ast
\bo_{\cR}
\bigg\|_{\infty}
.
\ee
The first 
factor
contributes $|\gW|^{2}$, and the second is $\ll|\gW|$, since $\gL_{0}$ is integrable.
For the last term, write
\bea
\nonumber
\bigg\|
\gL'_{k}\ast
\bo_{\cR}
\bigg\|_{\infty}
&\ll&
\sup_{z}
\int_{\R}
2^{2k}
|\gL'_{0}(u2^{2k})|
\bo_{\cR}
(z-u)
du
\\
\nonumber
&\ll&
\sum_{m\ge0}
2^{2(k-m)}
\sup_{|\cU|=2^{m+1-2k}}
|\cR\cap\cU|
,
\\
\label{eq:J2k1max}
\eea
where
the supremum is taken over intervals $\cU$.
Here we used that $\gL_{0}'$ has rapid decay, so certainly
$$
\gL_{0}'(y)\ll
\bo_{|y|<1}
+
\frac1{2^{2}}
\bo_{|y|<2}
+
\frac1{4^{2}}
\bo_{|y|<4}
+
\frac1{8^{2}}
\bo_{|y|<8}
+\cdots
.
$$

Note by \eqref{eq:apkBnd}
 that 
$$
|\cU|\ge2^{-2k+1}
\ge
\frac2{N^{1+\vep}}
.
$$ 

We have yet to exploit the structure of $\cR$ and do so now. This requires
the following
\begin{lem}\label{lem:ap2.3}
For any interval $\cU$ of length at least $1/N^{1+\vep}$, we have
\be\label{eq:ap2.3}
|\cR\cap\cU|
\ll
N^{-(1-\gd)+\vep}
|\cU|^{\gd+\vep}
.
\ee
\end{lem}

Postponing the proof of this lemma, we see that applying \eqref{eq:ap2.3} in \eqref{eq:J2k1max} gives
\be\label{eq:LinfgL}
\bigg\|
\gL'_{k}\ast
\bo_{\cR}
\bigg\|_{\infty}
\ll_{\vep}
2^{2k(1-\gd)}
N^{-(1-\gd)+\vep}
.
\ee
Putting \eqref{eq:LinfgL} into \eqref{eq:J2k1} gives \eqref{eq:ap3.8}, as claimed.
\epf

It remains to establish \eqref{eq:ap2.3}.

\pf[Proof of Lemma \ref{lem:ap2.3}]
From the structure of $\cR$ in \eqref{eq:ap2.1}, we have that
\beann
|\cR\cap\cU|
&\le&
\sum_{\ell}
\left|\cU\cap\bigg(\cI_{\ell}+[-2/N,2/N]\bigg)\right| %
\\
&\ll& 
\frac1N\#\left\{\ell\le L:\cU\cap\bigg(\cI_{\ell}+[-2/N,2/N]\bigg)\neq\emptyset\right\}
\\
&\ll_{\vep}& 
\frac1N
{
\mu(\cU+[0,1/N])
\over
N^{-\gd-\vep}
}
\\
&\ll_{\vep}& 
N^{-(1-\gd)+\vep}
|\cU|^{\gd+\vep}
,
\eeann
where we used \eqref{eq:ap1.2} and \eqref{eq:ap1.1} in the penultimate and final lines, respectively.
\epf

\subsubsection{The range $k> \cK$}\

In this range, we will establish
\begin{lem}
For $k> \cK$, and any $\vep>0$,
\be\label{eq:ap3.9}
\left|\cJ_{2}^{(k)}\right|
\ll_{\vep}
|\gW|^{3/2}\
2^{-k\gd}\
N^{-(1-\gd)+\vep}
.
\ee
\end{lem}


\pf
Changing variables $y\mapsto yx$ in \eqref{eq:ap3.5}, we have
$$
\left|
\cJ_{2}^{(k)}
\right|
\ll
\left|
\iint
\ff(x)
(\gL_{k}\ast\cS_{\gW})(xy)\
\eta(y)\
\bigg(
\gL'_{k}\ast
\bo_{\cR}
\bigg)
\left(y-1\right)\
 dy\
dx
\right|
,
$$
where we set
\be\label{eq:apFFis}
\ff(x):=x\ \eta(x)
(\bo_{\gW}\ast \bo_{-\gW})(-x)\
.
\ee

By the rapid decay of $\gL_{k}'$ and \eqref{eq:apRcont}, we may restrict the integral to $y\asymp1$ with a negligible error.
Now reverse orders, apply Parseval in 
$x$
,
reverse orders again, 
use the definition of the Fourier multipliers $\gl_{k}$, 
apply Cauchy-Schwarz in $y$, 
change variables $y\mapsto\xi/y$, 
and 
estimate:
\beann
\left|
\cJ_{2}^{(k)}
\right|
&\ll&
\left|
\int_{y\asymp1}
\left(
\int_{|\xi|/y\asymp2^{2k}}
\hat\ff(\xi)\
\overline{\hat\cS_{\gW}(\xi/y)}
\frac1y
d\xi
\right)
(\gL'_{k}\ast
\bo_{\cR})
\left(y-1\right)
dy
\right|
\\
&\ll&
\int_{|\xi|\asymp2^{2k}}
|\hat\ff(\xi)|\
\left(
\frac1{|\xi|}
\int
|\hat\cS_{\gW}(
y)|^{2}
dy
\right)^{1/2}
\|\gL'_{k}\ast
\bo_{\cR}\|_{2}\
d\xi
\\
&\ll_{\vep}&
\left(
2^{-k}
\int_{\xi\in\R}
|\hat\ff(\xi)|
d\xi
\right)
|\gW|
^{1/2}\
2^{k(1-\gd)}
N^{-(1-\gd)+\vep}
,
\eeann
where we estimated the last piece by
\beann
\left\|
\gL'_{k}\ast
\bo_{\cR}
\right\|_{2}
&\le&
\left\|
\gL'_{k}\ast
\bo_{\cR}
\right\|_{\infty}^{1/2}
\left\|
\gL'_{k}\ast
\bo_{\cR}
\right\|_{1}^{1/2}
\\
&\ll_{\vep}&
2^{k(1-\gd)}
N^{-(1-\gd)/2+\vep}
|\cR|^{1/2}
\\
&\ll_{\vep}&
2^{k(1-\gd)}
N^{-(1-\gd)+\vep}
,
\eeann
using the $\cL^{\infty}$ bound in \eqref{eq:LinfgL} and \eqref{eq:cRsize}. We easily estimate from \eqref{eq:apFFis} that $\|\hat\ff\|_{1}\ll|\gW|$, giving \eqref{eq:ap3.9}, as claimed.
\epf

\subsubsection{Completion of Proof}\

It is now a simple matter to establish Proposition \ref{prop:Jupper}. 
Putting \eqref{eq:ap3.8}, \eqref{eq:ap3.9}, and \eqref{eq:cJ2k} into \eqref{eq:J2J2p} gives
\bea
\nonumber
\cJ_{2}
&\ll_{\vep}&
\frac1M
N^{-(1-\gd)+\vep}
\left(
|\gW|^{3}\
2^{2\cK(1-\gd)}\
+
|\gW|^{3/2}\
2^{-\cK\gd}\
\right)
\\
\label{eq:cJ2bnd}
&\ll_{\vep}&
\frac1M
N^{-(1-\gd)+\vep}\
|\gW|^{3/(2-\gd)}\
,
\eea
on setting
\be\label{eq:cKIs}
\cK:=
{-3\log_{2}|\gW|\over 2(2-\gd)}
.
\ee
%
Combining \eqref{eq:cJ2bnd} with \eqref{eq:ap3.3} and choosing 
$$
M=
N^{-(1-\gd)/2+\vep}
|\gW|^{(1+\gd)/(4-2\gd)}
$$
gives \eqref{eq:Jupper}, as claimed. 
This competes the proof of Theorem \ref{thm:gdPlus}.

\bibliographystyle{alpha}

\bibliography{AKbibliog}

\begin{thebibliography}{MVW84}

\bibitem[BG08]{BourgainGamburd2008}
Jean Bourgain and Alex Gamburd.
\newblock Uniform expansion bounds for {C}ayley graphs of {${\rm SL}_2(\Bbb
  F_p)$}.
\newblock {\em Ann. of Math. (2)}, 167(2):625--642, 2008.

\bibitem[BGS10]{BourgainGamburdSarnak2010}
Jean Bourgain, Alex Gamburd, and Peter Sarnak.
\newblock Affine linear sieve, expanders, and sum-product.
\newblock {\em Invent. Math.}, 179(3):559--644, 2010.

\bibitem[BGS11]{BourgainGamburdSarnak2011}
J.~Bourgain, A.~Gamburd, and P.~Sarnak.
\newblock Generalization of {S}elberg's 3/16th theorem and affine sieve.
\newblock {\em Acta Math}, 207:255--290, 2011.

\bibitem[BK10]{BourgainKontorovich2010}
J.~Bourgain and A.~Kontorovich.
\newblock On representations of integers in thin subgroups of {SL}$(2,{{\bf
  {Z}}})$.
\newblock {\em GAFA}, 20(5):1144--1174, 2010.

\bibitem[BK11]{BourgainKontorovich2011}
J.~Bourgain and A.~Kontorovich.
\newblock On {Z}aremba's conjecture.
\newblock {\em Comptes Rendus Mathematique}, 349(9):493--495, 2011.

\bibitem[BKS10]{BourgainKontorovichSarnak2010}
J.~Bourgain, A.~Kontorovich, and P.~Sarnak.
\newblock Sector estimates for hyperbolic isometries.
\newblock {\em GAFA}, 20(5):1175--1200, 2010.

\bibitem[Bou03]{Bourgain2003}
J.~Bourgain.
\newblock On the {E}rd{\H o}s-{V}olkmann and {K}atz-{T}ao ring conjectures.
\newblock {\em Geom. Funct. Anal.}, 13(2):334--365, 2003.

\bibitem[Bou05]{Bourgain2005}
J.~Bourgain.
\newblock {\em Green's function estimates for lattice {S}chr\"odinger operators
  and applications}, volume 158 of {\em Annals of Mathematics Studies}.
\newblock Princeton University Press, Princeton, NJ, 2005.

\bibitem[Bou10]{Bourgain2010}
Jean Bourgain.
\newblock The discretized sum-product and projection theorems.
\newblock {\em Journal d'Analyse Math{\'e}matique}, 112:193--236, 2010.
\newblock 10.1007/s11854-010-0028-x.

\bibitem[Bum85]{Bumby1985}
Richard~T. Bumby.
\newblock Hausdorff dimension of sets arising in number theory.
\newblock In {\em Number theory ({N}ew {Y}ork, 1983--84)}, volume 1135 of {\em
  Lecture Notes in Math.}, pages 1--8. Springer, Berlin, 1985.

\bibitem[BV11]{BourgainVarju2011}
J.~Bourgain and P.~Varj{\'u}.
\newblock Expansion in ${SL}\sb n(\bold {Z}/q \bold {Z})$, $q$ arbitrary, 2011.
\newblock To appear, {\it Invent. Math.} {\tt arXiv:1006.3365v1}.

\bibitem[Dol98]{Dolgopyat1998}
Dmitry Dolgopyat.
\newblock On decay of correlations in {A}nosov flows.
\newblock {\em Ann. of Math. (2)}, 147(2):357--390, 1998.

\bibitem[Goo41]{Good1941}
I.~J. Good.
\newblock The fractional dimensional theory of continued fractions.
\newblock {\em Proc. Cambridge Philos. Soc.}, 37:199--228, 1941.

\bibitem[Goo83]{GoodBook}
Anton Good.
\newblock {\em Local analysis of {S}elberg's trace formula}, volume 1040 of
  {\em Lecture Notes in Mathematics}.
\newblock Springer-Verlag, Berlin, 1983.

\bibitem[GS01]{GoldsteinSchlag2001}
Michael Goldstein and Wilhelm Schlag.
\newblock H\"older continuity of the integrated density of states for
  quasi-periodic {S}chr\"odinger equations and averages of shifts of
  subharmonic functions.
\newblock {\em Ann. of Math. (2)}, 154(1):155--203, 2001.

\bibitem[GV11]{GolsefidyVarju2011}
A.~Golsefidy and P.~Varj{\'u}, 2011.
\newblock Preprint.

\bibitem[Haa81]{Haagerup1981}
Uffe Haagerup.
\newblock The best constants in the {K}hintchine inequality.
\newblock {\em Studia Math.}, 70(3):231--283 (1982), 1981.

\bibitem[Hen89]{Hensley1989}
Doug Hensley.
\newblock The distribution of badly approximable numbers and continuants with
  bounded digits.
\newblock In {\em Th\'eorie des nombres ({Q}uebec, {PQ}, 1987)}, pages
  371--385. de Gruyter, Berlin, 1989.

\bibitem[Hen92]{Hensley1992}
Doug Hensley.
\newblock Continued fraction {C}antor sets, {H}ausdorff dimension, and
  functional analysis.
\newblock {\em J. Number Theory}, 40(3):336--358, 1992.

\bibitem[Hen96]{Hensley1996}
Douglas Hensley.
\newblock A polynomial time algorithm for the {H}ausdorff dimension of
  continued fraction {C}antor sets.
\newblock {\em J. Number Theory}, 58(1):9--45, 1996.

\bibitem[Hen06]{Hensley2006}
Doug Hensley.
\newblock {\em Continued fractions}.
\newblock World Scientific Publishing Co. Pte. Ltd., Hackensack, NJ, 2006.

\bibitem[Jen04]{Jenkinson2004}
Oliver Jenkinson.
\newblock On the density of {H}ausdorff dimensions of bounded type continued
  fraction sets: the {T}exan conjecture.
\newblock {\em Stoch. Dyn.}, 4(1):63--76, 2004.

\bibitem[JP01]{JenkinsonPollicott2001}
Oliver Jenkinson and Mark Pollicott.
\newblock Computing the dimension of dynamically defined sets: {$E_2$} and
  bounded continued fractions.
\newblock {\em Ergodic Theory Dynam. Systems}, 21(5):1429--1445, 2001.

\bibitem[Kon13]{Kontorovich2013}
Alex Kontorovich.
\newblock From {A}pollonius to {Z}aremba: local-global phenomena in thin
  orbits.
\newblock {\em Bull. Amer. Math. Soc. (N.S.)}, 50(2):187--228, 2013.

\bibitem[KS03]{KimSarnak2003}
H.~Kim and P.~Sarnak.
\newblock Refined estimates towards the {R}amanujan and {S}elberg conjectures.
\newblock {\em J. Amer. Math. Soc.}, 16(1):175--181, 2003.

\bibitem[Lal89]{Lalley1989}
Steven~P. Lalley.
\newblock Renewal theorems in symbolic dynamics, with applications to geodesic
  flows, non-{E}uclidean tessellations and their fractal limits.
\newblock {\em Acta Math.}, 163(1-2):1--55, 1989.

\bibitem[McM09]{McMullen2009}
Curtis~T. McMullen.
\newblock Uniformly {D}iophantine numbers in a fixed real quadratic field.
\newblock {\em Compos. Math.}, 145(4):827--844, 2009.

\bibitem[MVW84]{MatthewsVasersteinWeisfeiler1984}
C.~Matthews, L.~Vaserstein, and B.~Weisfeiler.
\newblock Congruence properties of {Z}ariski-dense subgroups.
\newblock {\em Proc. London Math. Soc}, 48:514--532, 1984.

\bibitem[Nau05]{Naud2005}
Fr{\'e}d{\'e}ric Naud.
\newblock Expanding maps on {C}antor sets and analytic continuation of zeta
  functions.
\newblock {\em Ann. Sci. \'Ecole Norm. Sup. (4)}, 38(1):116--153, 2005.

\bibitem[Nie78]{Niederreiter1978}
Harald Niederreiter.
\newblock Quasi-{M}onte {C}arlo methods and pseudo-random numbers.
\newblock {\em Bull. Amer. Math. Soc.}, 84(6):957--1041, 1978.

\bibitem[Zar72]{Zaremba1972}
S.~K. Zaremba.
\newblock La m\'ethode des ``bons treillis'' pour le calcul des int\'egrales
  multiples.
\newblock In {\em Applications of number theory to numerical analysis ({P}roc.
  {S}ympos., {U}niv. {M}ontreal, {M}ontreal, {Q}ue., 1971)}, pages 39--119.
  Academic Press, New York, 1972.

\end{thebibliography}

\end{document}